\documentclass[11pt, reqno]{amsart} \textwidth=14.9cm \oddsidemargin=1cm \evensidemargin=1cm
\usepackage{amsmath,amssymb,amscd,mathrsfs,stmaryrd,wasysym}
\usepackage{amsmath}
\usepackage{amsxtra}
\usepackage{amscd}
\usepackage{amsthm}
\usepackage{amsfonts}
\usepackage{amssymb}
\usepackage{eucal}
\usepackage{color}
\usepackage{graphicx}%
\newtheorem{theorem}{Theorem}[section]
\newtheorem{lemma}[theorem]{Lemma}

\theoremstyle{definition}
\newtheorem{definition}[theorem]{Definition}
\newtheorem{remark}[theorem]{Remark}

\definecolor{A}{rgb}{.75,1,.75}
\def\nslash{\:\notslash\:}

\numberwithin{equation}{section}

\begin{document}

\title[A categorical equivalence]{A categorical equivalence between affine Yokonuma-\\Hecke algebras and some quiver Hecke algebras}
\author[Weideng Cui]{Weideng Cui}
\address{School of Mathematics, Shandong University, Jinan, Shandong 250100, P.R. China.}
\email{cwdeng@amss.ac.cn}

\begin{abstract}
Inspired by the work of Rostam, we establish an explicit categorical equivalence between affine Yokonuma-Hecke algebras and quiver Hecke algebras associated to disjoint copies of quivers of (affine) type $A,$ generalizing Rouquier's categorical equivalence theorem.
\end{abstract}



\maketitle
\medskip
\section{Introduction}
\subsection{}
The representation theory of affine and degenerate affine Hecke algebras is very important, and has been studied extensively over the past few decades. Kazhdan and Lusztig \cite{KaLu} and Ginzburg \cite{CG} gave a classification and construction of irreducible representations of the affine Hecke algebra $H_{q,k}$ when $k$ is the complex number field and $q$ is not a root of unity, which is known as the Deligne-Langlands-Lusztig classification. When $H_{q,k}$ is of affine type $A,$ a classification of irreducible representations of $H_{q,k}$ was obtained in \cite{AM} for any $q$ and arbitrary sufficiently large $k.$ When the parameter $q$ is not a root of the Poincar\'{e} polynomial, Xi \cite{Xi} proved that the classification established in \cite{KaLu} remains valid. Lusztig \cite{Lu1} proved that certain module category of an affine Hecke algebra is equivalent to its suitable counterpart of some graded affine Hecke algebra.

\subsection{}
In order to compute the decomposition numbers of the Iwahori-Hecke algebra of type $A$ at an $e$th root of unity over $\mathbb{C}$, Lascoux, Leclerc, and Thibon \cite{LLT} suggested the conjecture that the decomposition numbers can be computed using the canonical bases of quantum groups of affine type $A.$ Later on, using the geometric results of Kazhdan-Lusztig and Ginzburg, Ariki \cite{Ari} proved this conjecture by establishing a connection between the representation theory of cyclotomic Hecke algebras and the canonical bases of integrable highest weight modules over $\widehat{\mathfrak{sl}}_{e}(\mathbb{C}).$ Ariki's work turned out to be a beginning of some exciting developments which continue to this day, namely categorification.

In order to provide a categorification of quantum groups, Khovanov and Lauda \cite{KhLa1, KhLa2} and Rouquier \cite{Rou1} have independently introduced a new family of algebras associated to a generalized Cartan matrix, known as quiver Hecke algebras or Khovanov-Lauda-Rouquier algebras. Moreover, they have shown that there exists an algebra isomorphism between the integral form of the negative half of a quantum group and the Grothendieck group of the additive category of finitely generated graded projective modules of a quiver Hecke algebra. Later on, when the Cartan matrix is symmetric, Varagnolo and Vasserot \cite{VV} and independently Rouquier \cite{Rou2} proved that Lusztig's canonical bases or Kashiwara's lower global bases correspond to the isomorphism classes of graded self-dual indecomposable projective modules under this isomorphism.

Besides, Khovanov and Lauda \cite{KhLa1} suggested the cyclotomic categorification conjecture, which was proved by Brundan and Kleshchev \cite{BK2} for type $A_{\infty}$ and $A_{n}^{(1)}$ using the isomorphism between the cyclotomic quiver Hecke algebras of type $A$ and the cyclotomic Hecke algebras of type $A$ which was established in \cite{BK1}, and proved by Kang and Kashiwara \cite{KK} for all symmetrizable Kac-Moody algebras. Brundan and Kleshchev's isomorphism theorem was independently proved by Rouquier [Rou1, Corollaries 3.17 and 3.20], which allows us to construct an explicit $\mathbb{Z}$-grading on the cyclotomic Hecke algebras of type $A$ and study their graded representation theory \cite{BK2}. Rouquier \cite{Rou1} also proved that certain module category of an affine and degenerate affine Hecke algebra is equivalent to its suitable counterpart of some quiver Hecke algebra of (affine) type $A$.


\subsection{}
Yokonuma-Hecke algebras were introduced by Yokonuma \cite{Yo} as a centralizer algebra associated to the permutation representation of a finite Chevalley group $G$ with respect to a maximal unipotent subgroup of $G$. The Yokonuma-Hecke algebra $Y_{r,n}(q)$ (of type $A$) is a quotient of the group algebra of the modular framed braid group $(\mathbb{Z}/r\mathbb{Z})\wr B_{n},$ where $B_{n}$ is the braid group on $n$ strands (of type $A$). By the presentation given by Juyumaya and Kannan \cite{Ju1, Ju2, JuK}, the Yokonuma-Hecke algebra $Y_{r,n}(q)$ can also be regraded as a deformation of the group algebra of the complex reflection group $G(r,1,n),$ which is isomorphic to the wreath product $(\mathbb{Z}/r\mathbb{Z})\wr \mathfrak{S}_{n}$. It is well-known that there exists another deformation of the group algebra of $G(r,1,n),$ namely the Ariki-Koike algebra \cite{AK}. The Yokonuma-Hecke algebra $Y_{r,n}(q)$ is quite different from the Ariki-Koike algebra.  For example, the Iwahori-Hecke algebra of type $A$ is canonically a subalgebra of the Ariki-Koike algebra, whereas it is an obvious quotient of $Y_{r,n}(q),$ but not an obvious subalgebra of it.

Recently, by generalizing the approach of Okounkov-Vershik \cite{OV} on the representation theory of $\mathfrak{S}_n$, Chlouveraki and Poulain d'Andecy \cite{ChP1} introduced the notion of affine Yokonuma-Hecke algebra $\widehat{Y}_{r,n}(q)$ and gave explicit formulas for all irreducible representations of $Y_{r,n}(q)$ over $\mathbb{C}(q)$, and obtained a semisimplicity criterion for it. In their subsequent paper [ChP2], they studied the representation theory of the affine Yokonuma-Hecke algebra $\widehat{Y}_{r,n}(q)$ and the cyclotomic Yokonuma-Hecke algebra $Y_{r,n}^{d}(q)$. In particular, they gave the classification of irreducible representations of $Y_{r,n}^{d}(q)$ in the generic semisimple case. In \cite{CW}, we gave the classification of the simple $\widehat{Y}_{r,n}(q)$-modules as well as the classification of the simple modules of the cyclotomic Yokonuma-Hecke algebras over an algebraically closed field $\mathbb{K}$ of characteristic $p$ such that $p$ does not divide $r.$ In the past several years, the study of affine and cyclotomic Yokonuma-Hecke algebras has made substantial progress; see \cite{ChP1, ChP2, ChS, C, CW, ER, JP, Lu2, Ro}.

\subsection{} Recently, Rostam \cite{Ro} has shown that cyclotomic Yokonuma-Hecke algebras are particular cases of cyclotomic quiver Hecke algebras, generalizing the results of Brundan and Kleshchev. Largely inspired by the work of Rostam, we establish an explicit categorical equivalence between affine Yokonuma-Hecke algebras $\widehat{Y}_{r,n}(q)$ and quiver Hecke algebras $H_{n}(\Gamma)$ associated to disjoint copies of quivers of (affine) type $A$, thus generalizing Rouquier's categorical equivalence theorem [Rou1, Theorems 3.16 and 3.19].

This paper is organized as follows. In Section 2, we establish many necessary results and then state the main theorem \ref{Morita-equi}. In Section 3, we give the proof of the main theorem by verifying the defining relations of $\widehat{Y}_{r,n}(q)$ and $H_{n}(\Gamma),$ respectively. In Section 4, we consider the degenerate case.

\section{An equivalence of module categories}

\subsection{Quiver Hecke algebras}
Let $r, n\in \mathbb{Z}_{\geq 1},$ and let $k$ be a field which contains a primitive $r$th root of unity $\zeta.$ Let $p$ be the characteristic of $k$ and consider an element $q\in k^{*}.$ We denote by $e$ the smallest integer such that $1+q+\cdots+q^{e-1}=0$, and set $e:=\infty$ if no such integer exists.

Assume first that $q=1.$ Given a subset $I$ of $k,$ we denote by $I_{1}$ the quiver with set of vertices $I$ and with an arrow $i\rightarrow i+1$ whenever $i, i+1\in I.$

Assume now that $q\neq 1.$ Given a subset $I$ of $k^{*},$ we denote by $I_{q}$ the quiver with set of vertices $I$ and with an arrow $i\rightarrow qi$ whenever $i, qi\in I.$

Assume that $I_{q}$ is connected. Note that when $q=1,$ $I_{q}$ has type $A_{\ell}$ if $\vert I\vert=\ell< p$ or $|I|=\ell$ and $p=0;$ type $A_{\infty}$ if $|I|=\infty$ and $p=0;$ type $\widetilde{A}_{p-1}$ if $|I|=p> 0.$ When $q\neq 1,$ $I_{q}$ has type $A_{\ell}$ if $\vert I\vert=\ell< e$ or $|I|=\ell$ and $e=\infty;$ type $A_{\infty}$ if $|I|=\infty$ and $e=\infty;$ type $\widetilde{A}_{e-1}$ if $|I|=e> 0.$

Let $J=\{0,1,\ldots,r-1\}$. We now define the quiver $\Gamma=\coprod_{j\in J} I_{q},$ which is $r$ disjoint copies of $I_{q}.$ Hence, the vertices of $\Gamma$ are the elements of $K=I\times J$ and there is an arrow between $(i, j)$ and $(i', j')$ in $\Gamma$ if and only if there is an arrow between $i$ and $i'$ in $I_{q}$ and $j=j'.$ For $k, k'\in K,$ we write $k\rightleftarrows k'$ whenever $k\rightarrow k'$ and $k'\rightarrow k$; this happens only when $e=2.$ Finally, we write $k\nslash k'$ if $k\neq k'$ and if there is no arrow between these two vertices. The action of $\mathfrak{S}_{n}$ on $K^{n}=I^{n}\times J^{n}$ is given by the diagonal action, that is, $\sigma\cdot(\nu, j)=(\sigma(\nu), \sigma(j))$ for any $(\nu, j)\in K^{n}$ and $\sigma\in \mathfrak{S}_{n}.$

\begin{definition}
The quiver Hecke algebra $H_{n}(\Gamma)$ associated to $\Gamma$ is a (possibly non-unitary) $k$-algebra with generators $\{1_{\alpha}\}_{\alpha\in K^{n}}$, $\{x_{a}\}_{1\leq a\leq n}$ and $\{\tau_{i}\}_{1\leq i\leq n-1}$ and relations:

\begin{align}
1_{\alpha}1_{\alpha'}&=\delta_{\alpha,\alpha'}1_{\alpha},\qquad\sum_{\alpha\in K^{n}}1_{\alpha}=1;\label{rel-def-QHA1}\\[0.2em]
\tau_{i}1_{\alpha}&=1_{s_{i}(\alpha)}\tau_{i},\qquad x_{a}1_{\alpha}=1_{\alpha}x_{a};\label{rel-def-QHA2}\\[0.2em]
x_{a}x_{b}&=x_{b}x_{a},\qquad \quad \tau_{i}\tau_{j}=\tau_{j}\tau_{i} \quad \mbox{if $\vert i-j\vert \geq 2$;}\label{rel-def-QHA3}\\[0.2em]
\tau_{i}x_{a}1_{\alpha}-x_{s_{i}(a)}\tau_{i}1_{\alpha}&=\left\{
\begin{array}{ll}
-1_{\alpha}&\text{if $a=i$ and $\alpha_i = \alpha_{i+1}$},\\[0.2em]
1_{\alpha}&\text{if $a=i+1$ and $\alpha_i = \alpha_{i+1}$},\\[0.2em]
0&\text{otherwise;}
\end{array}
\right.\label{rel-def-QHA4}\\[0.2em]
\tau_{i}^{2}1_{\alpha}&=Q_{\alpha_i, \alpha_{i+1}}(x_{i}, x_{i+1})1_{\alpha};\label{rel-def-QHA5}\\[0.2em]
\nonumber\tau_{i+1}\tau_{i}\tau_{i+1}1_{\alpha}-\tau_{i}\tau_{i+1}\tau_{i}1_{\alpha}&=\end{align}
\begin{align}
\left\{
\begin{array}{ll}
(x_{i+2}-x_{i})^{-1}(Q_{\alpha_i, \alpha_{i+1}}(x_{i+2}, x_{i+1})-Q_{\alpha_i, \alpha_{i+1}}(x_{i}, x_{i+1}))1_{\alpha}&\text{if $\alpha_i = \alpha_{i+2}$},\\[0.3em]
0&\text{otherwise,}
\end{array}
\right.\label{rel-def-QHA6}
\end{align}
for $\alpha,\alpha'\in K^{n},$ $1\leq a,b\leq n$ and $1\leq i,j\leq n-1,$ and where\\
\[Q_{\alpha_i, \alpha_{i+1}}(u, v)=\left\{
\begin{array}{ll}
0&\text{if $\alpha_i = \alpha_{i+1}$},\\[0.05em]
1&\text{if $\alpha_i\nslash \alpha_{i+1}$},\\[0.05em]
v-u&\text{if $\alpha_i\rightarrow \alpha_{i+1}$},\\[0.05em]
u-v&\text{if $\alpha_i\leftarrow \alpha_{i+1}$},\\[0.05em]
-(u-v)^{2}&\text{if $\alpha_i\rightleftarrows \alpha_{i+1}$}.
\end{array}
\right.\]
\end{definition}

By \eqref{rel-def-QHA4}, we get the following useful identity in $H_{n}(\Gamma)$:\vspace{0.15cm}
\begin{align}
f\tau_{i}1_{(\nu, j)}=\left\{
\begin{array}{ll}
\tau_{i}{}^{s_{i}}f1_{(\nu, j)}+\partial_{i}(f)1_{(\nu, j)}&\text{if $(\nu_i, j_i)=(\nu_{i+1}, j_{i+1})$},\\[0.2em]
\tau_{i}{}^{s_{i}}f1_{(\nu, j)}&\text{otherwise},
\end{array}
\right.\label{commu-relations}
\end{align}
\vspace{0.15cm}
for any $1\leq i\leq n-1,$ $(\nu, j)\in K^{n}$ and $f\in k[[x_1,\dots,x_n]]$, and where
\begin{equation*}
\partial_{i}(f):=\frac{f-{}^{s_{i}}f}{x_{i+1}-x_{i}}.
\end{equation*}

In the rest of this section, we always assume that $q\neq 1.$ For each $(\nu, j)\in K^{n}$ and $1\leq i\leq n-1,$ we define the following element:
\[Q_{i}(\nu, j)=
\begin{cases}
\begin{cases}
q(x_{i+1}+1)-(x_{i}+1) & \text{if } \nu_i = \nu_{i+1},
\\
\left(q^{-1}(x_{i}+1)-(x_{i+1}+1)\right)^{-1} & \text{if } \nu_{i+1} =q\nu_{i},
\\
\frac{\nu_{i+1}(x_{i+1}+1)-q\nu_{i}(x_{i}+1)}
{\nu_{i}(x_{i}+1)-\nu_{i+1}(x_{i+1}+1)}
 & \text{otherwise,}
\end{cases}
& \text{if } j_i = j_{i+1},
\\
f_{i, j} & \text{if } j_i \neq j_{i+1},
\end{cases}
\]
\setlength{\parskip}{2pt}where $f_{i,j}\in \{1, q\}$ is given by

\[f_{i, j}=
\begin{cases}
q & \text{if } j_i < j_{i+1},
\\[0.1em]
1& \text{if } j_i > j_{i+1},
\end{cases}
\]
and

\[P_{i}(\nu, j)=
\begin{cases}
\begin{cases}
1 & \text{if } \nu_i = \nu_{i+1},
\\
\frac{(q-1)\nu_{i+1}(x_{i+1}+1)}{\nu_{i}(x_{i}+1)-\nu_{i+1}(x_{i+1}+1)}
 & \text{otherwise,}
\end{cases}
& \text{if } j_i = j_{i+1},
\\
0 & \text{if } j_i \neq j_{i+1}.
\end{cases}
\]

The following lemma can be easily checked by definition.
\begin{lemma}
\begin{align} P_{i}(\nu, j)+{}^{s_{i}}P_{i}(s_{i}\cdot(\nu, j))&=1-q \quad\text{ if }\nu_i \neq \nu_{i+1} \text{ and } j_i = j_{i+1};\label{P-relations1}\\[0.2em]
{}^{s_{i}}P_{i+1}(s_{i}\cdot(\nu, j))&={}^{s_{i+1}}P_{i}(s_{i+1}\cdot(\nu, j))\quad\text{ for any }(\nu, j)\in K^{n}. \label{P-relations2}
\end{align}
\end{lemma}
The following lemma can also be easily checked by definition.
\begin{lemma}
\begin{align}
{}^{s_{i}}Q_{i+1}(s_{i+1}s_{i}\cdot(\nu, j))&={}^{s_{i+1}}Q_{i}(s_{i}s_{i+1}\cdot(\nu, j))\quad\text{ for any }(\nu, j)\in K^{n}. \label{Q-relations1}
\end{align}
\end{lemma}
We also need the following key lemma.

\begin{lemma}
\begin{align}
\hspace{-0.1cm}{}^{s_{i}}Q_{i}(s_{i}\cdot(\nu, j))Q_{i}(\nu, j)=\left\{
\begin{array}{ll}
\hspace{-0.15cm}(1-P_{i}(\nu, j))(q+P_{i}(\nu, j))&\text{if $(\nu_i, j_{i})\nslash (\nu_{i+1}, j_{i+1})$},\\[0.3em]
\hspace{-0.15cm}\frac{(1-P_{i}(\nu, j))(q+P_{i}(\nu, j))}{x_{i+1}-x_{i}}&\text{if $(\nu_i, j_{i})\rightarrow (\nu_{i+1}, j_{i+1})$},\\[0.3em]
\hspace{-0.15cm}\frac{(1-P_{i}(\nu, j))(q+P_{i}(\nu, j))}{x_{i}-x_{i+1}}&\text{if $(\nu_i, j_{i})\leftarrow (\nu_{i+1}, j_{i+1})$},\\[0.3em]
\hspace{-0.15cm}\frac{(1-P_{i}(\nu, j))(q+P_{i}(\nu, j))}{(x_{i}-x_{i+1})(x_{i+1}-x_{i})}&\text{if $(\nu_i, j_{i})\rightleftarrows (\nu_{i+1}, j_{i+1})$}.\label{Q-relations2}
\end{array}
\right.\end{align}
\end{lemma}

\begin{proof}
(1) If $j_{i}\neq j_{i+1}$, then $(\nu_i, j_{i})\nslash (\nu_{i+1}, j_{i+1})$, in this case \[{}^{s_{i}}Q_{i}(s_{i}\cdot(\nu, j))Q_{i}(\nu, j)=f_{i,j}f_{i, s_{i}(j)}=q,\] while $(1-P_{i}(\nu, j))(q+P_{i}(\nu, j))=q.$\vspace{0.25cm}

(2) If $j_{i}= j_{i+1}$ and $(\nu_i, j_{i})\nslash (\nu_{i+1}, j_{i+1}),$ then we have
\begin{align*}
{}^{s_{i}}Q_{i}(s_{i}\cdot(\nu, j))Q_{i}(\nu, j)&=\frac{\nu_{i}(x_{i}+1)-q\nu_{i+1}(x_{i+1}+1)}
{\nu_{i+1}(x_{i+1}+1)-\nu_{i}(x_{i}+1)}\cdot\frac{\nu_{i+1}(x_{i+1}+1)-q\nu_{i}(x_{i}+1)}
{\nu_{i}(x_{i}+1)-\nu_{i+1}(x_{i+1}+1)}\\
&=\frac{\big(\nu_{i}(x_{i}+1)-q\nu_{i+1}(x_{i+1}+1)\big)\cdot\big(q\nu_{i}(x_{i}+1)-\nu_{i+1}(x_{i+1}+1)\big)}{\big(\nu_{i}(x_{i}+1)-\nu_{i+1}(x_{i+1}+1)\big)^{2}},
\end{align*}
while
\begin{align}
&(1-P_{i}(\nu, j))(q+P_{i}(\nu, j))\notag\\=&q-\frac{(q-1)^{2}\nu_{i+1}(x_{i+1}+1)}{\nu_{i}(x_{i}+1)-\nu_{i+1}(x_{i+1}+1)}-
\frac{(q-1)^{2}\nu_{i+1}^{2}(x_{i+1}+1)^{2}}{\big(\nu_{i}(x_{i}+1)-\nu_{i+1}(x_{i+1}+1)\big)^{2}}\notag\\
=&\big\{q\big(\nu_{i}(x_{i}+1)-\nu_{i+1}(x_{i+1}+1)\big)^{2}-(q-1)^{2}\nu_{i+1}(x_{i+1}+1)\big(\nu_{i}(x_{i}+1)-
\nu_{i+1}(x_{i+1}+1)\big)\notag\\
&-(q-1)^{2}\nu_{i+1}^{2}(x_{i+1}+1)^{2}\big\}\cdot\big(\nu_{i}(x_{i}+1)-\nu_{i+1}(x_{i+1}+1)\big)^{-2}\notag\\
=&\frac{q\nu_{i}^{2}(x_{i}+1)^{2}+q\nu_{i+1}^{2}(x_{i+1}+1)^{2}-(q^{2}+1)\nu_{i}\nu_{i+1}(x_{i}+1)(x_{i+1}+1)}
{\big(\nu_{i}(x_{i}+1)-\nu_{i+1}(x_{i+1}+1)\big)^{2}}\notag\\
=&\frac{\big(\nu_{i}(x_{i}+1)-q\nu_{i+1}(x_{i+1}+1)\big)\cdot\big(q\nu_{i}(x_{i}+1)-\nu_{i+1}(x_{i+1}+1)\big)}
{\big(\nu_{i}(x_{i}+1)-\nu_{i+1}(x_{i+1}+1)\big)^{2}}.\label{Q-relations3}
\end{align}\vspace{0.25cm}

(3) If $j_{i}= j_{i+1}$ and $(\nu_i, j_{i})\rightarrow (\nu_{i+1}, j_{i+1}),$ then $\nu_{i+1}=q\nu_{i},$ and we have
\begin{align*}
{}^{s_{i}}Q_{i}(s_{i}\cdot(\nu, j))Q_{i}(\nu, j)&=\frac{\nu_{i}(x_{i}+1)-q\nu_{i+1}(x_{i+1}+1)}
{\nu_{i+1}(x_{i+1}+1)-\nu_{i}(x_{i}+1)}\cdot\frac{1}{q^{-1}(x_{i}+1)-(x_{i+1}+1)}\\
&=\frac{q(x_{i+1}+1)-q^{-1}(x_{i}+1)}{\big(q^{-1}(x_{i}+1)-(x_{i+1}+1)\big)^{2}},
\end{align*}
by setting $\nu_{i+1}=q\nu_{i}$ in \eqref{Q-relations3}, we get
\begin{align*}
(1-P_{i}(\nu, j))(q+P_{i}(\nu, j))&=\frac{\big(\nu_{i}(x_{i}+1)-q^{2}\nu_{i}(x_{i+1}+1)\big)\cdot\big(q\nu_{i}(x_{i}+1)-q\nu_{i}(x_{i+1}+1)\big)}
{\big(\nu_{i}(x_{i}+1)-q\nu_{i}(x_{i+1}+1)\big)^{2}}\\
&=\frac{\big(q(x_{i+1}+1)-q^{-1}(x_{i}+1)\big)\cdot(x_{i+1}-x_{i})}{\big(q^{-1}(x_{i}+1)-(x_{i+1}+1)\big)^{2}}.
\end{align*}\vspace{0.25cm}

(4) If $j_{i}= j_{i+1}$ and $(\nu_i, j_{i})\leftarrow (\nu_{i+1}, j_{i+1}),$ then $\nu_{i}=q\nu_{i+1},$ and we have
\begin{align*}
{}^{s_{i}}Q_{i}(s_{i}\cdot(\nu, j))Q_{i}(\nu, j)&=\frac{1}{q^{-1}(x_{i+1}+1)-(x_{i}+1)}\cdot
\frac{\nu_{i+1}(x_{i+1}+1)-q\nu_{i}(x_{i}+1)}{\nu_{i}(x_{i}+1)-\nu_{i+1}(x_{i+1}+1)}\\
&=\frac{q(x_{i}+1)-q^{-1}(x_{i+1}+1)}{\big(q^{-1}(x_{i+1}+1)-(x_{i}+1)\big)^{2}},
\end{align*}
by setting $\nu_{i}=q\nu_{i+1}$ in \eqref{Q-relations3}, we get
\begin{align*}
(1-P_{i}(\nu, j))(q+P_{i}(\nu, j))&=\frac{\big(\nu_{i}(x_{i}+1)-\nu_{i}(x_{i+1}+1)\big)\cdot\big(q\nu_{i}(x_{i}+1)-q^{-1}\nu_{i}(x_{i+1}+1)\big)}
{\big(\nu_{i}(x_{i}+1)-q^{-1}\nu_{i}(x_{i+1}+1)\big)^{2}}\\
&=\frac{\big(q(x_{i}+1)-q^{-1}(x_{i+1}+1)\big)\cdot(x_{i}-x_{i+1})}{\big(q^{-1}(x_{i+1}+1)-(x_{i}+1)\big)^{2}}.
\end{align*}\vspace{0.25cm}

(5) If $j_{i}= j_{i+1}$ and $(\nu_i, j_{i})\rightleftarrows (\nu_{i+1}, j_{i+1}),$ then $q=-1$ and $\nu_{i+1}=-\nu_{i},$ and we have
\begin{align*}
{}^{s_{i}}Q_{i}(s_{i}\cdot(\nu, j))Q_{i}(\nu, j)&=\frac{1}{-(x_{i+1}+1)-(x_{i}+1)}\cdot
\frac{1}{-(x_{i}+1)-(x_{i+1}+1)}\\
&=\frac{1}{(x_{i}+x_{i+1}+2)^{2}},
\end{align*}
by setting $\nu_{i+1}=-\nu_{i}$ in \eqref{Q-relations3}, we get
\begin{align*}
(1-P_{i}(\nu, j))(q+P_{i}(\nu, j))&=\frac{\nu_{i}(x_{i}-x_{i+1})\cdot\nu_{i}(x_{i+1}-x_{i})}
{\big(\nu_{i}(x_{i}+x_{i+1}+2)\big)^{2}}\\
&=\frac{(x_{i}-x_{i+1})\cdot(x_{i+1}-x_{i})}{(x_{i}+x_{i+1}+2)^{2}}.
\end{align*}
We are done.\end{proof}

By \eqref{rel-def-QHA5} and \eqref{Q-relations2}, we get the following consequence: for any $(\nu, j)\in K^{n},$
\begin{equation}
(\tau_{i}^{2}){}^{s_{i}}Q_{i}(s_{i}\cdot(\nu, j))Q_{i}(\nu, j)1_{(\nu, j)}=(1-P_{i}(\nu, j))(q+P_{i}(\nu, j))1_{(\nu, j)}.\label{Q-P-relations}
\end{equation}

\subsection{Affine Yokonuma-Hecke algebras}
Let $q$ be an indeterminate and let $\mathcal{R}=\mathbb{Z}[\frac{1}{r}][q,q^{-1}].$ We consider $k$ as an $\mathcal{R}$-algebra by mapping $q$ to the invertible element $q\in k^{*}.$

\begin{definition}
The affine Yokonuma-Hecke algebra, denoted by $\widehat{Y}_{r,n}(q)$, is an $\mathcal{R}$-associative algebra generated by the elements $t_{1},\ldots,t_{n},g_{1},\ldots,g_{n-1},X_{1}^{\pm1},$ in which the generators $t_{1},\ldots,t_{n},g_{1},$ $\ldots,g_{n-1}$ satisfy the following relations:
\begin{align}
g_ig_j&=g_jg_i \qquad \qquad\qquad\quad\hspace{0.3cm}\mbox{for all $i,j=1,\ldots,n-1$ such that $\vert i-j\vert \geq 2$;}\label{rel-def-Y1}\\[0.1em]
g_ig_{i+1}g_i&=g_{i+1}g_ig_{i+1} \qquad \quad\qquad\hspace{0.05cm}\mbox{for all $i=1,\ldots,n-2$;}\label{rel-def-Y2}\\[0.1em]
t_it_j&=t_jt_i \qquad\qquad\qquad\qquad  \mbox{for all $i,j=1,\ldots,n$;}\label{rel-def-Y3}\\[0.1em]
g_it_j&=t_{s_i(j)}g_i \quad \quad\qquad\qquad\hspace{0.25cm}\mbox{for all $i=1,\ldots,n-1$ and $j=1,\ldots,n$;}\label{rel-def-Y4}\\[0.1em]
t_i^r&=1 \quad \qquad\qquad\qquad\qquad\mbox{for all $i=1,\ldots,n$;}\label{rel-def-Y5}\\[0.2em]
g_{i}^{2}&=q+(q-1)e_{i}g_{i} \qquad \hspace{0.48cm}\mbox{for all $i=1,\ldots,n-1$;}\label{rel-def-Y6}
\end{align}
where $s_{i}$ is the transposition $(i,i+1)$, and for each $1\leq i\leq n-1$,
$$e_{i} :=\frac{1}{r}\sum\limits_{s=0}^{r-1}t_{i}^{s}t_{i+1}^{-s},$$
together with the following relations concerning the generators $X_{1}^{\pm1}$:
\begin{align}
X_{1}X_{1}^{-1}&=X_{1}^{-1}X_{1}=1;\label{rel-def-Y7}\\[0.1em]
g_{1}X_{1}g_{1}X_{1}&=X_{1}g_{1}X_{1}g_{1};\label{rel-def-Y8}\\[0.1em]
g_{i}X_{1}&=X_{1}g_{i} \qquad \qquad \quad\mbox{for all $i=2,\ldots,n-1$;}\label{rel-def-Y9}\\[0.1em]
t_{j}X_{1}&=X_{1}t_{j} \qquad \qquad\quad\mbox{for all $j=1,\ldots,n$;}\label{rel-def-Y10}
\end{align}
\end{definition}

We define inductively elements $X_{2},\ldots,X_{n}$ in $\widehat{Y}_{r,n}(q)$ by
\begin{equation}
X_{i+1} :=q^{-1}g_{i}X_{i}g_{i}\quad\mathrm{for}~i=1,\ldots,n-1.\label{X2n}
\end{equation}
Then it is proved in [ChP1, Lemma 1] that we have, for any $1\leq i\leq n-1$,
\begin{equation}
g_{i}X_{j}=X_{j}g_{i}\quad\mathrm{for}~j=1,2,\ldots,n~\mathrm{such~that}~j\neq i, i+1.\label{giXj}
\end{equation}
Moreover, by [ChP1, Proposition 1], we have that the elements $t_{1},\ldots, t_{n}, X_{1},\ldots, X_{n}$ form a commutative family, that is,
\begin{equation}
xy=yx\quad\mathrm{for~any}~x,y\in \{t_{1},\ldots, t_{n}, X_{1},\ldots, X_{n}\}.\label{xyyx}
\end{equation}
We also have the following identities (see [ChP2, Lemma 2.3]): for $1\leq i\leq n-1$,
\begin{align}
g_{i}X_{i+1}&=X_{i}g_{i}+(q-1)e_{i}X_{i+1},\label{gxxg}\\
X_{i+1}g_{i}&=g_{i}X_{i}+(q-1)e_{i}X_{i+1},\label{gxxg1}\\
g_{i}X_{i}^{-1}&=X_{i+1}^{-1}g_{i}+(q-1)e_{i}X_{i}^{-1}, \label{gxxg2}\\
X_{i}^{-1}g_{i}&=g_{i}X_{i+1}^{-1}+(q-1)e_{i}X_{i}^{-1}.\label{gxxg3}
\end{align}

Let $\widehat{Y}_{r,n}^{k}=k\otimes_{\mathcal{R}}\widehat{Y}_{r,n}(q)$ and let $M$ be a $\widehat{Y}_{r,n}^{k}$-module. Given $\alpha=(\nu, j)\in K^{n},$ we denote by $M_{\alpha}$ the subspace of $M$ on which $X_{a}-\nu_{a}$ acts locally nilpotently for $1\leq a\leq n,$ and simultaneously, $t_{a}-\zeta^{j_{a}}$ acts as zero for $1\leq a\leq n.$ Let $\mathcal{C}_{K}$ be the category of finitely generated $\widehat{Y}_{r,n}^{k}$-modules $M$ such that $M=\bigoplus_{\alpha\in K^{n}}M_{\alpha}.$

Given an object $M\in \mathcal{C}_{K}$, we can write $M$ as the direct sum of its weight spaces (simultaneous generalized eigenspaces):
\begin{equation}
M(\nu, j)=\{v\in M\:|\:(X_{a}-\nu_{a})^{N}v=(t_{a}-\zeta^{j_{a}})v=0\text{ for all }1\leq a\leq n\text{ and }N\gg 0\}.\label{M-J}
\end{equation}
Considering the weight space decomposition of the regular module, we deduce that there is a family $\{e(\nu,j)\}_{(\nu,j)\in K^{n}}$ of mutually orthogonal idempotents in $\widehat{Y}_{r,n}^{k}$ such that $e(\nu,j)M=M(\nu, j)$ for each $M\in \mathcal{C}_{K}$. In fact, each $e(\nu, j)$ lies in the commutative subalgebra generated by $X_{1}^{\pm1},\ldots, X_{n}^{\pm1}, t_{1},\ldots,t_{n}$, all but finitely many of the $e(\nu, j)$'s are zero, and their sum is the identity element in $\widehat{Y}_{r,n}^{k}.$

By definition, we easily get the following equalities: for $1\leq i\leq n-1$ and $(\nu, j)\in K^{n}$,
\begin{equation}
e_{i}e(\nu, j)=0 \text{ if } j_{i}\neq j_{i+1} ~~~~\text{ and }~~~~e_{i}e(\nu, j)=e(\nu, j) \text{ if }j_{i}=j_{i+1}.\label{eee}
\end{equation}

From \eqref{gxxg}, \eqref{gxxg1} and \eqref{eee}, we immediately get the following lemma.

\begin{lemma}
For $1\leq i\leq n-1$ and $(\nu, j)\in K^{n}$ with $j_{i}\neq j_{i+1},$ we have
\begin{equation}
g_{i}X_{i+1}e(\nu, j)=X_{i}g_{i}e(\nu, j) ~~~~\text{ and }~~~~X_{i+1}g_{i}e(\nu, j)=g_{i}X_{i}e(\nu, j).\label{gxe-xge}
\end{equation}
\end{lemma}

We also have the following lemma.
\begin{lemma}
If $1\leq i\leq n-1$ and $(\nu, j)\in K^{n}$ with $j_{i}\neq j_{i+1},$ then we have
\begin{equation}
g_{i}e(\nu, j)=e(s_{i}\cdot(\nu, j))g_{i}.\label{ge-eg}
\end{equation}
\end{lemma}

\begin{proof}
From \eqref{rel-def-Y4} and \eqref{gxe-xge}, we see that $g_{i}e(\nu, j)$ maps $M(\nu, j)$ to $M(s_{i}\cdot(\nu, j))$ for any $M\in \mathcal{C}_{K}.$
Thus, both $g_{i}e(\nu, j)$ and $e(s_{i}\cdot(\nu, j))g_{i}$ map $M(\nu', j')$ to zero unless $(\nu', j')=(\nu, j)$, and they map each $m\in M(\nu, j)$ to $g_{i}m.$ We get \eqref{ge-eg}.
\end{proof}

For $1\leq i\leq n-1,$ set
\[\Phi_i=g_i + (1-q)\sum_{\substack{(\nu, j) \in K^{n} \\ \nu_i \neq \nu_{i+1} \\ j_i = j_{i+1}}} \left(1 - X_i X_{i+1}^{-1}\right)^{-1} e(\nu, j) + \sum_{\substack{(\nu, j) \in K^{n} \\ \nu_i = \nu_{i+1} \\ j_i = j_{i+1}}} e(\nu, j).
\]

Recall also that the following intertwining element introduced in [CW, $\S5.1$]:
$$\Theta_{i} :=g_{i}(1-X_{i}X_{i+1}^{-1})+(1-q)e_{i},~~~~1\leq i\leq n-1.$$

\begin{lemma}
For $1\leq i, k\leq n-1$ and $(\nu, j)\in K^{n}$, we have
\begin{equation}
\Theta_{i}\Theta_{i+1}\Theta_{i}=\Theta_{i+1}\Theta_{i}\Theta_{i+1},~~~~~~\Theta_{i}\Theta_{k}=\Theta_{k}\Theta_{i}\text{ if } \vert i-k\vert> 1;\label{theta}
\end{equation}
\begin{equation}
\Phi_{i}e(\nu, j)=e(s_{i}\cdot(\nu, j))\Phi_{i}.\label{phie-ephi}
\end{equation}
\end{lemma}

\begin{proof}
The second identity in \eqref{theta} can be easily proved. The first identity in \eqref{theta} can be proved by making a lengthy but routine calculation. Here we provide a sketch of another proof. To prove it, it suffices to prove that $\Theta_{i}\Theta_{i+1}\Theta_{i}e(\nu, j)=\Theta_{i+1}\Theta_{i}\Theta_{i+1}e(\nu, j)$ for all $(\nu, j)\in K^{n},$ which can be proved by discussing the following five cases:\vspace{0.2cm}

Case 1: $j_i=j_{i+1}=j_{i+2};$

Case 2: $j_i=j_{i+1}\neq j_{i+2};$

Case 3: $j_i\neq j_{i+1}=j_{i+2};$

Case 4: $j_i=j_{i+2}\neq j_{i+1};$

Case 5: $j_i, j_{i+1}, j_{i+2}$ are all different.\vspace{0.2cm}

By $e_{i}e(\nu, j)=e_{i+1}e(\nu, j)=e(\nu, j),$ Case 1 comes down to the case of affine Hecke algebra of type $A,$ which is exactly [BK, (4.13)].

For Case 2, after a careful calculation, we get that
\begin{align*}\Theta_{i}\Theta_{i+1}\Theta_{i}e(\nu, j)=&g_{i}g_{i+1}g_{i}(1-X_{i+1}X_{i+2}^{-1})(1-X_{i}X_{i+2}^{-1})(1-X_{i}X_{i+1}^{-1})e(\nu, j)\\&+(1-q)g_{i}(1-X_{i}X_{i+1}^{-1})g_{i+1}(1-X_{i+1}X_{i+2}^{-1})e(\nu, j),\end{align*}
and
\begin{align*}\Theta_{i+1}\Theta_{i}\Theta_{i+1}e(\nu, j)=&g_{i+1}g_{i}g_{i+1}(1-X_{i}X_{i+1}^{-1})(1-X_{i}X_{i+2}^{-1})(1-X_{i+1}X_{i+2}^{-1})e(\nu, j)\\&+(1-q)g_{i}(1-X_{i}X_{i+1}^{-1})g_{i+1}(1-X_{i+1}X_{i+2}^{-1})e(\nu, j).\end{align*}
Thus, Case 2 holds by the braid relation \eqref{rel-def-Y2}.

Case 3 holds similarly.

For Case 4, after a careful calculation, we get that
\begin{align*}\Theta_{i}\Theta_{i+1}\Theta_{i}e(\nu, j)=&g_{i}g_{i+1}g_{i}(1-X_{i+1}X_{i+2}^{-1})(1-X_{i}X_{i+2}^{-1})(1-X_{i}X_{i+1}^{-1})e(\nu, j)\\&+(1-q)g_{i}^{2}(1-X_{i+1}X_{i+2}^{-1})(1-X_{i}X_{i+1}^{-1})e(\nu, j),\end{align*}
and
\begin{align*}\Theta_{i+1}\Theta_{i}\Theta_{i+1}e(\nu, j)=&g_{i+1}g_{i}g_{i+1}(1-X_{i}X_{i+1}^{-1})(1-X_{i}X_{i+2}^{-1})(1-X_{i+1}X_{i+2}^{-1})e(\nu, j)\\&+(1-q)g_{i+1}^{2}(1-X_{i}X_{i+1}^{-1})(1-X_{i+1}X_{i+2}^{-1})e(\nu, j).\end{align*}
Thus, Case 4 holds by the braid relation \eqref{rel-def-Y2} and the fact that $g_{i}^{2}e(\nu, j)=g_{i+1}^{2}e(\nu, j)=qe(\nu, j)$.

For Case 5, after a careful calculation, we get that
\begin{align*}\Theta_{i}\Theta_{i+1}\Theta_{i}e(\nu, j)=g_{i}g_{i+1}g_{i}(1-X_{i+1}X_{i+2}^{-1})(1-X_{i}X_{i+2}^{-1})(1-X_{i}X_{i+1}^{-1})e(\nu, j),\end{align*}
and
\begin{align*}\Theta_{i+1}\Theta_{i}\Theta_{i+1}e(\nu, j)=g_{i+1}g_{i}g_{i+1}(1-X_{i}X_{i+1}^{-1})(1-X_{i}X_{i+2}^{-1})(1-X_{i+1}X_{i+2}^{-1})e(\nu, j).\end{align*}
Thus, Case 5 holds by the braid relation \eqref{rel-def-Y2}.

We now prove \eqref{phie-ephi}. If $j_{i}=j_{i+1},$ the identity follows from [BK, (4.14)]; if $j_{i}\neq j_{i+1},$ the identity follows from \eqref{ge-eg}.
\end{proof}

We also need the following lemma.

\begin{lemma}
For each $(\nu, j)\in K^{n}$, we have
\begin{equation}
\Phi_{r}X_{s}=X_{s}\Phi_{r}\text{ if } s\neq r, r+1;\label{phi-x}
\end{equation}
\begin{equation}
\Phi_{r}\Phi_{s}=\Phi_{s}\Phi_{r}\text{ if } \vert r-s\vert> 1;\label{phi-phi}
\end{equation}
\begin{equation}
\Phi_{r}Q_{s}'(\nu, j)=Q_{s}'(\nu, j)\Phi_{r}\text{ if } \vert r-s\vert> 1;\label{phi-q}
\end{equation}
\begin{equation}
\Phi_{r}X_{r+1}e(\nu, j)=X_{r}\Phi_{r}e(\nu, j)\text{ if } \nu_{r}\neq \nu_{r+1}\text{ and } j_{r}=j_{r+1};\label{phi-x-nuj1}
\end{equation}
\begin{equation}
X_{r+1}\Phi_{r}e(\nu, j)=\Phi_{r}X_{r}e(\nu, j)\text{ if } \nu_{r}\neq \nu_{r+1}\text{ and } j_{r}=j_{r+1};\label{phi-x-nuj2}
\end{equation}
\begin{equation}
\Phi_r^2e(\nu, j)=
\frac{(X_{r+1}-qX_r)(X_r-qX_{r+1})}{(X_{r+1}-X_{r})(X_{r}-X_{r+1})}e(\nu, j)\text{ if } \nu_{r}\neq \nu_{r+1}\text{ and } j_{r}=j_{r+1};\label{phi-2-nuj1}
\end{equation}
\begin{align*}
\hspace{-3cm}\Phi_r\Phi_{r+1}\Phi_re(\nu, j)=&\end{align*}
\begin{align}\left\{
\begin{array}{ll}
\hspace{-0.15cm}(\Phi_{r+1}\Phi_r\Phi_{r+1}+q\Phi_{r}-q\Phi_{r+1})e(\nu, j)&\hspace{-0.3cm}\text{ if } \nu_r=\nu_{r+2}=\nu_{r+1},\\[0.2em]
\hspace{-0.15cm}(\Phi_{r+1}\Phi_r\Phi_{r+1}+Z_r)e(\nu, j)&\hspace{-0.3cm}\text{ if } \nu_r=\nu_{r+2}\neq \nu_{r+1},~~~~\text{if } j_r=j_{r+1}=j_{r+2},\\[0.2em]\label{phi-braid}
\hspace{-0.15cm}\Phi_{r+1}\Phi_r\Phi_{r+1}e(\nu, j)&\hspace{-0.15cm}\text{otherwise,}
\end{array}
\right.
\end{align}
where $Z_r$ denotes $(1-q)^2\frac{(X_rX_{r+2}-X_{r+1}^2)(X_r X_{r+1}-qX_{r+1}X_{r+2})}{(X_r-X_{r+1})^2(X_{r+1}-X_{r+2})^2}.$
\end{lemma}

\begin{proof}
The identities \eqref{phi-x} and \eqref{phi-phi} are clear from definitions. The identity \eqref{phi-q} follows from \eqref{phi-x}. The identities \eqref{phi-x-nuj1} and \eqref{phi-x-nuj2} follow from [CW, (5.2)] and the fact that $\Phi_r e(\nu, j) = \Theta_r (1-X_{r}X_{r+1}^{-1})^{-1} e(\nu, j)$ for $\nu_{r}\neq \nu_{r+1}$ and $j_{r}=j_{r+1}$, or they follow from [BK, (4.17) and (4.18)]. The identities \eqref{phi-2-nuj1} and \eqref{phi-braid} follow from [BK, (4.19) and (4.20)].\end{proof}

\subsection{An equivalence of categories}
Recall that $\mathcal{C}_{K}$ is the category of finitely generated $\widehat{Y}_{r,n}^{k}$-modules $M$ such that $M=\bigoplus_{\alpha\in K^{n}}M_{\alpha}.$ We denote by $\mathcal{C}_{\Gamma}^{0}$ the category of finitely generated $H_{n}(\Gamma)$-modules $M$ such that for every $\alpha\in K^{n},$ $x_{a}1_{\alpha}$ acts locally nilpotently on $1_{\alpha}M$ for $1\leq a\leq n.$

Now we can state the main theorem of this paper.
\begin{theorem}\label{Morita-equi}
There is an equivalence of categories $\Phi: \mathcal{C}_{\Gamma}^{0}\rightarrow \mathcal{C}_{K}$ given by $M\mapsto M,$ and where for each $\alpha=(\nu, j)\in K^{n},$ $t_{a}$ acts on $1_{\alpha}M$ by $\zeta^{j_{a}},$ $X_{a}$ acts on $1_{\alpha}M$ by $\nu_{a}(x_{a}+1)$ for $1\leq a\leq n,$ and $g_{i}$ acts on $1_{\alpha}M$ by $\tau_{i}Q_{i}(\nu, j)-P_{i}(\nu, j)$ for $1\leq i\leq n-1.$

The inverse functor $\Psi: \mathcal{C}_{K}\rightarrow \mathcal{C}_{\Gamma}^{0}$ is given by $N\mapsto N,$ where for each $\alpha=(\nu, j)\in K^{n},$ $1_{\alpha}$ acts on $N$ by $e(\nu, j),$ $x_{a}$ acts on $N_{\alpha}$ by $\nu_{a}^{-1}X_{a}-1$ for $1\leq a\leq n,$ and $\tau_{i}$ acts on $N_{\alpha}$ by $\Phi_{i}Q_{i}'(\nu, j)^{-1}=(g_{i}+P_{i}'(\nu, j))Q_{i}'(\nu, j)^{-1}$ for $1\leq i\leq n-1,$ where $P_{i}'(\nu, j)$ and $Q_{i}'(\nu, j)$ are defined by replacing $x_{i}$ with $\nu_{i}^{-1}X_{i}-1$ in the expressions of $P_{i}(\nu, j)$ and $Q_{i}(\nu, j).$
\end{theorem}

\section{Proof of the main theorem}

In this section we prove the main theorem \ref{Morita-equi}, which is motivated by the work of [BK1] and [Ro].\vspace{0.3cm}

{\it Proof of Theorem \ref{Morita-equi}} We first verify that the actions under $\Phi$ of the generators of $\widehat{Y}_{r,n}^{k}$ satisfy the relations \eqref{rel-def-Y1}-\eqref{rel-def-Y10}.

The relation \eqref{rel-def-Y1} easily follows from \eqref{rel-def-QHA3} and \eqref{rel-def-QHA4}.

The relations \eqref{rel-def-Y3} and \eqref{rel-def-Y5} easily follow from the definition.

The relation \eqref{rel-def-Y4} easily follows from the next fact:
\begin{equation*}
\zeta^{j_{b}}P_{a}(\nu, j)=\zeta^{j_{s_{a}(b)}}P_{a}(\nu, j)\quad\text{ for any }1\leq a\leq n-1, ~1\leq b\leq n\text{ and }(\nu, j)\in K^{n}.
\end{equation*}

Next we check the relation \eqref{rel-def-Y6}, that is, $g_{i}^{2}\diamond 1_{(\nu, j)}m=(q+(q-1)g_{i}e_{i})\diamond 1_{(\nu, j)}m$ for all $(\nu, j)\in K^{n}$ and $m\in M$ with an $M\in \mathcal{C}_{\Gamma}^{0}.$ To simplify notation for the remainder of the proof, we no longer write the elements ``$m$" on the right hand side of all expressions, but remember it is always there. By definition, we have
\begin{align*}
&g_{i}^{2}\diamond 1_{(\nu, j)}\\=&(\tau_{i}Q_{i}(s_{i}\cdot(\nu, j))\tau_{i}Q_{i}(\nu, j)-P_{i}(s_{i}\cdot(\nu, j))\tau_{i}Q_{i}(\nu, j)-\tau_{i}Q_{i}(\nu, j)P_{i}(\nu, j)+P_{i}(\nu, j)^{2})1_{(\nu, j)}.
\end{align*}

If $(\nu_{i}, j_{i})=(\nu_{i+1}, j_{i+1}),$ by $\tau_{i}^{2}1_{(\nu, j)}=0,$ $e_{i}\diamond 1_{(\nu, j)}=1$ and $\partial_{i}(qx_{i+1}+q-x_{i}-1)=q+1$ and \eqref{commu-relations}, we deduce that
\begin{align*}
g_{i}^{2}\diamond 1_{(\nu, j)}&=(\tau_{i}(qx_{i+1}+q-x_{i}-1)\tau_{i}Q_{i}(\nu, j)-2\tau_{i}Q_{i}(\nu, j)+1)1_{(\nu, j)}\\
&=((q+1)\tau_{i}Q_{i}(\nu, j)-2\tau_{i}Q_{i}(\nu, j)+1)1_{(\nu, j)}\\
&=(q+(q-1)g_{i}e_{i})\diamond 1_{(\nu, j)}.
\end{align*}

If $j_{i}=j_{i+1}$ and $\nu_{i}\neq \nu_{i+1},$ by $e_{i}\diamond 1_{(\nu, j)}=1$, \eqref{commu-relations}, \eqref{P-relations1} and \eqref{Q-P-relations}, we have
\begin{align*}
&g_{i}^{2}\diamond 1_{(\nu, j)}\\=&((\tau_{i}^{2}){}^{s_{i}}Q_{i}(s_{i}\cdot(\nu, j))Q_{i}(\nu, j)-\tau_{i}(P_{i}(\nu, j)+{}^{s_{i}}P_{i}(s_{i}\cdot(\nu, j)))Q_{i}(\nu, j)+P_{i}(\nu, j)^{2})1_{(\nu, j)}\\
=&((1-P_{i}(\nu, j))(q+P_{i}(\nu, j))+(q-1)\tau_{i}Q_{i}(\nu, j)+P_{i}(\nu, j)^{2})1_{(\nu, j)}\\
=&(q+(q-1)g_{i}e_{i})\diamond 1_{(\nu, j)}.
\end{align*}

If $j_{i}\neq j_{i+1},$ by $e_{i}\diamond 1_{(\nu, j)}=0$, $P_{i}(\nu, j)=0$, \eqref{commu-relations} and \eqref{Q-P-relations}, we have
\begin{align*}
g_{i}^{2}\diamond 1_{(\nu, j)}&=(\tau_{i}^{2}){}^{s_{i}}Q_{i}(s_{i}\cdot(\nu, j))Q_{i}(\nu, j)1_{(\nu, j)}\\
&=(1-P_{i}(\nu, j))(q+P_{i}(\nu, j))1_{(\nu, j)}\\
&=q1_{(\nu, j)}=(q+(q-1)g_{i}e_{i})\diamond 1_{(\nu, j)}.
\end{align*}

To prove \eqref{rel-def-Y8}, it suffices to show that \eqref{X2n} holds. Since we have proved \eqref{rel-def-Y6}, it suffices to show that \eqref{gxxg} holds.

If $(\nu_{i}, j_{i})=(\nu_{i+1}, j_{i+1}),$ by $e_{i}\diamond 1_{(\nu, j)}=1$ and \eqref{commu-relations}, we deduce that
\begin{align*}
X_{i}g_{i}\diamond 1_{(\nu, j)}&=X_{i}\diamond (\tau_{i}Q_{i}(\nu, j)-P_{i}(\nu, j))1_{(\nu, j)}\\
&=\nu_{i+1}(x_{i}+1)\tau_{i}Q_{i}(\nu, j)1_{(\nu, j)}-\nu_{i}(x_{i}+1)P_{i}(\nu, j)1_{(\nu, j)}\\
&=(\nu_{i}\tau_{i}(x_{i+1}+1)-\nu_{i})Q_{i}(\nu, j)1_{(\nu, j)}-\nu_{i}(x_{i}+1)P_{i}(\nu, j)1_{(\nu, j)}\\
&=\nu_{i}\tau_{i}(x_{i+1}+1)Q_{i}(\nu, j)1_{(\nu, j)}-\nu_{i}(q(x_{i+1}+1)-(x_{i}+1)+(x_{i}+1))1_{(\nu, j)}\\
&=\nu_{i}\tau_{i}(x_{i+1}+1)Q_{i}(\nu, j)1_{(\nu, j)}-q\nu_{i}(x_{i+1}+1)1_{(\nu, j)},
\end{align*}
and
\begin{align*}
\hspace{-0.6cm}(g_{i}+(1-q)e_{i})X_{i+1}\diamond 1_{(\nu, j)}&=(\tau_{i}Q_{i}(\nu, j)-P_{i}(\nu, j)+(1-q))\nu_{i+1}(x_{i+1}+1)1_{(\nu, j)}\\
&=\nu_{i}\tau_{i}(x_{i+1}+1)Q_{i}(\nu, j)1_{(\nu, j)}-q\nu_{i}(x_{i+1}+1)1_{(\nu, j)}.
\end{align*}

If $j_{i}=j_{i+1}$ and $\nu_{i}\neq \nu_{i+1},$ by $e_{i}\diamond 1_{(\nu, j)}=1$ and \eqref{commu-relations}, we have
\begin{align*}
X_{i}g_{i}\diamond 1_{(\nu, j)}&=\nu_{i+1}(x_{i}+1)\tau_{i}Q_{i}(\nu, j)1_{(\nu, j)}-\nu_{i}(x_{i}+1)P_{i}(\nu, j)1_{(\nu, j)}\\
&=\nu_{i+1}\tau_{i}(x_{i+1}+1)Q_{i}(\nu, j)1_{(\nu, j)}-\nu_{i}(x_{i}+1)P_{i}(\nu, j)1_{(\nu, j)},
\end{align*}
and
\begin{align*}
&(g_{i}+(1-q)e_{i})X_{i+1}\diamond 1_{(\nu, j)}\\=&(\tau_{i}Q_{i}(\nu, j)-P_{i}(\nu, j)+(1-q))\nu_{i+1}(x_{i+1}+1)1_{(\nu, j)}\\
=&\nu_{i+1}\tau_{i}(x_{i+1}+1)Q_{i}(\nu, j)1_{(\nu, j)}+\nu_{i+1}(x_{i+1}+1)(1-q)\big(1+\frac{\nu_{i+1}(x_{i+1}+1)}{\nu_{i}(x_{i}+1)-\nu_{i+1}(x_{i+1}+1)}\big)\\
=&\nu_{i+1}\tau_{i}(x_{i+1}+1)Q_{i}(\nu, j)1_{(\nu, j)}-\nu_{i}(x_{i}+1)P_{i}(\nu, j)1_{(\nu, j)}.
\end{align*}

If $j_{i}\neq j_{i+1},$ by $e_{i}\diamond 1_{(\nu, j)}=0$, $P_{i}(\nu, j)=0$ and \eqref{commu-relations}, we have
\begin{align*}
X_{i}g_{i}\diamond 1_{(\nu, j)}=\nu_{i+1}\tau_{i}(x_{i+1}+1)Q_{i}(\nu, j)1_{(\nu, j)},
\end{align*}
and
\begin{align*}
(g_{i}+(1-q)e_{i})X_{i+1}\diamond 1_{(\nu, j)}=\nu_{i+1}\tau_{i}(x_{i+1}+1)Q_{i}(\nu, j)1_{(\nu, j)}.
\end{align*}

To prove \eqref{rel-def-Y9}, it suffices to prove \eqref{giXj}, which easily follows from \eqref{rel-def-QHA4}.

To prove \eqref{rel-def-Y10}, it suffices to prove \eqref{xyyx}, which immediately follows from the definition and \eqref{rel-def-QHA2}.

Finally we need to check the braid relations (2.2). Without loss of generality we assume that $i=1$ and $n=3,$ and we need to show that $g_{2}g_{1}g_{2}\diamond 1_{(\nu, j)}=g_{1}g_{2}g_{1}\diamond 1_{(\nu, j)},$ where $\nu=(\nu_1, \nu_2, \nu_3)$, $j=(j_1, j_2, j_3)$ and $(\nu, j)\in K^{3}.$

If $j_1=j_2=j_3$, we set $r:=\nu_1, s:=\nu_2, t:=\nu_3.$ We stop writing $\diamond 1_{(\nu, j)}$ on the right hand side of all expressions and stop writing $j=(j_1, j_2, j_3)$ in the expressions of $Q_{a}(\nu, j)$ and $P_{a}(\nu, j)$, that is, we write $Q_{a}(\nu)$ instead of $Q_{a}(\nu, j)$ and $P_{a}(\nu)$ instead of $P_{a}(\nu, j)$ for $a=1,2.$ By definition, $g_{2}g_{1}g_{2}$ and $g_{1}g_{2}g_{1}$ equal
\begin{equation}\label{ELPQ}
\begin{split}
&-P_2(rst)P_1(rst)P_2(rst)+ \tau_2Q_2(rst)P_1(rst)P_2(rst)
\\
&+ P_2(srt)\tau_1Q_1(rst)P_2(rst)
-\tau_2Q_2(srt)\tau_1Q_1(rst)P_2(rst)
\\
&+ P_2(rts)P_1(rts)\tau_2Q_2(rst)
- \tau_2Q_2(rts)P_1(rts)\tau_2Q_2(rst)
\\
&- P_2(trs)\tau_1Q_1(rts)\tau_2Q_2(rst)
+ \tau_2Q_2(trs)\tau_1Q_1(rts)\tau_2Q_2(rst),
\end{split}
\end{equation}
and
\begin{equation}\label{ERPQ}
\begin{split}
&-P_1(rst)P_2(rst)P_1(rst)+ \tau_1Q_1(rst)P_2(rst)P_1(rst)
\\&+ P_1(rts)\tau_2Q_2(rst)P_1(rst)
-\tau_1Q_1(rts)\tau_2Q_2(rst)P_1(rst)
\\&+ P_1(srt)P_2(srt)\tau_1Q_1(rst)
- \tau_1Q_1(srt)P_2(srt)\tau_1Q_1(rst)
\\& - P_1(str)\tau_2Q_2(srt)\tau_1Q_1(rst)
+ \tau_1Q_1(str)\tau_2Q_2(srt)\tau_1Q_1(rst),
\end{split}
\end{equation}
respectively. We have to show that \eqref{ELPQ} equals \eqref{ERPQ}. In order to show this, we compare various $\tau$-coefficients after commuting all $\tau$'s to the left using \eqref{commu-relations}. Then we need to consider five cases.

\subsubsection*{Case 1: $r, s, t$ are all different.}

\

\vspace{0.35cm}\noindent By \eqref{commu-relations}, \eqref{ELPQ} equals

\begin{equation*}
\begin{split}
&-P_2(rst)P_1(rst)P_2(rst)+ \tau_2Q_2(rst)P_1(rst)P_2(rst)
\\
&+ \tau_1{}^{s_{1}}P_2(srt)Q_1(rst)P_2(rst)
-\tau_2\tau_1{}^{s_{1}}Q_2(srt)Q_1(rst)P_2(rst)
\\
&+ \tau_2{}^{s_{2}}P_2(rts){}^{s_{2}}P_1(rts)Q_2(rst)
- (\tau_2^{2}){}^{s_{2}}Q_2(rts){}^{s_{2}}P_1(rts)Q_2(rst)
\\
&- \tau_1\tau_2{}^{s_{2}s_{1}}P_2(trs){}^{s_{2}}Q_1(rts)Q_2(rst)
+ \tau_2\tau_1\tau_2{}^{s_{2}s_{1}}Q_2(trs){}^{s_{2}}Q_1(rts)Q_2(rst),
\end{split}
\end{equation*}
and \eqref{ERPQ} equals
\begin{equation*}
\begin{split}
&-P_1(rst)P_2(rst)P_1(rst)+ \tau_1Q_1(rst)P_2(rst)P_1(rst)
\\&+ \tau_2{}^{s_{2}}P_1(rts)Q_2(rst)P_1(rst)
-\tau_1\tau_2{}^{s_{2}}Q_1(rts)Q_2(rst)P_1(rst)
\\&+ \tau_1{}^{s_{1}}P_1(srt){}^{s_{1}}P_2(srt)Q_1(rst)
- (\tau_1^{2}){}^{s_{1}}Q_1(srt){}^{s_{1}}P_2(srt)Q_1(rst)
\\& - \tau_2\tau_1{}^{s_{1}s_{2}}P_1(str){}^{s_{1}}Q_2(srt)Q_1(rst)
+ \tau_1\tau_2\tau_1{}^{s_{1}s_{2}}Q_1(str){}^{s_{1}}Q_2(srt)Q_1(rst).
\end{split}
\end{equation*}

Since $r,s, t$ are all different, we have $\tau_2\tau_1\tau_2=\tau_1\tau_2\tau_1,$ and their corresponding coefficients are equal to each other by \eqref{Q-relations1}. The corresponding $\tau_1\tau_2$-coefficients and $\tau_2\tau_1$-coefficients are equal to each other by \eqref{P-relations2}. For the $\tau_1$-coefficients, it suffices to prove that
\begin{equation*}
{}^{s_{1}}P_2(srt)P_2(rst)-P_2(rst)P_1(rst)={}^{s_{1}}P_1(srt){}^{s_{1}}P_2(srt).
\end{equation*}
We have
\begin{align*}
&{}^{s_{1}}P_2(srt)P_2(rst)-P_2(rst)P_1(rst)\\=&\big(\frac{t(x_{3}+1)}{r(x_{1}+1)-t(x_{3}+1)}-\frac{s(x_{2}+1)}{r(x_{1}+1)-s(x_{2}+1)}\big)\cdot \frac{t(x_{3}+1)}{s(x_{2}+1)-t(x_{3}+1)}\\
=&\frac{-rt(x_1+1)(x_3+1)}{(r(x_{1}+1)-t(x_{3}+1))(r(x_{1}+1)-s(x_{2}+1))},
\end{align*}
and
\begin{align*}
{}^{s_{1}}P_1(srt){}^{s_{1}}P_2(srt)&=\frac{r(x_1+1)}{s(x_{2}+1)-r(x_{1}+1)}\cdot \frac{t(x_{3}+1)}{r(x_{1}+1)-t(x_{3}+1)}\\
&=\frac{-rt(x_1+1)(x_3+1)}{(r(x_{1}+1)-t(x_{3}+1))(r(x_{1}+1)-s(x_{2}+1))}.
\end{align*}
The $\tau_2$-coefficients can be handed similarly. Considering the constant coefficients, we need to show that
\begin{align*}
P_1(rst)&P_2(rst)^{2}+(\tau_2^{2}){}^{s_{2}}Q_2(rts)Q_2(rst){}^{s_{2}}P_1(rts)\\
&=P_2(rst)P_1(rst)^{2}+(\tau_1^{2}){}^{s_{1}}Q_1(srt)Q_1(rst){}^{s_{1}}P_2(srt).
\end{align*}
Using the following identities:
$$
{}^{s_{2}}P_1(rts)={}^{s_{1}}P_2(srt),
$$$$
(\tau_2^{2}){}^{s_{2}}Q_2(rts)Q_2(rst)=(1-P_{2}(rst))(q+P_{2}(rst)),
$$
$$
(\tau_1^{2}){}^{s_{1}}Q_1(srt)Q_1(rst)=(1-P_{1}(rst))(q+P_{1}(rst)),
$$
we need to check that
\begin{align*}
P_2(rst)^{2}&(P_1(rst)-{}^{s_{2}}P_1(rts))+(1-q)P_2(rst){}^{s_{2}}P_1(rts)\\&=P_1(rst)^{2}(P_2(rst)-{}^{s_{1}}P_2(srt))+(1-q)P_1(rst){}^{s_{1}}P_2(srt).
\end{align*}
By an explicit expansion, it is easy to see that the left hand side and right hand side are both equal to
\begin{align*}
(1-q)^{3}\frac{st^{2}(x_2+1)(x_{3}+1)^{2}}{(t(x_{3}+1)-s(x_{2}+1))(t(x_{3}+1)-r(x_{1}+1))(s(x_{2}+1)-r(x_{1}+1))}.
\end{align*}

\subsubsection*{Case 2: $r=s\neq t$.}
\subsubsection*{Case 3: $r\neq s=t$.}
\subsubsection*{Case 4: $r= t\neq s$.}

\

\vspace{0.3cm}\noindent As in the proof of \emph{Case 1}, these three cases can be proved in exactly the same way as in the proof of [BK1, Theorem 4.3, p. 478-479]. We omit the details.

\subsubsection*{Case 5: $r=s=t$.}

\

\vspace{0.3cm}\noindent This case is left as an exercise in the proof of [BK1, Theorem 4.3, p. 479]. Here we shall give a brief proof. In this case, we have $P_{1}(\nu)=P_{2}(\nu)=1$, $Q_{1}(\nu)=qx_2+q-x_1-1$ and $Q_{2}(\nu)=qx_3+q-x_2-1.$ Moreover, By \eqref{commu-relations}, \eqref{ELPQ} equals
\begin{align*}
\begin{split}
&-1+\tau_2Q_2(\nu)+\tau_1Q_1(\nu)-\tau_2Q_2(\nu)\tau_1Q_1(\nu)+\tau_2Q_2(\nu)-\tau_2Q_2(\nu)\tau_2Q_2(\nu)\\&-\tau_1Q_1(\nu)\tau_2Q_2(\nu)
+\tau_2Q_2(\nu)\tau_1Q_1(\nu)\tau_2Q_2(\nu)\\=&-1+\tau_2Q_2(\nu)+\tau_1Q_1(\nu)-\tau_2\tau_1{}^{s_{1}}Q_2(\nu)Q_1(\nu)
-\tau_2\partial_{1}(Q_2(\nu))Q_1(\nu)+\tau_2Q_2(\nu)\\&-\tau_2\partial_{2}(Q_2(\nu))Q_2(\nu)-\tau_1\tau_2{}^{s_{2}}Q_1(\nu)Q_2(\nu)
-\tau_1\partial_{2}(Q_1(\nu))Q_2(\nu)\\&+\tau_2\tau_1\tau_2{}^{s_{2}s_{1}}Q_{2}(\nu){}^{s_{2}}Q_1(\nu)Q_2(\nu)+
\tau_2\tau_1\partial_{2}({}^{s_{1}}Q_2(\nu)){}^{s_{2}}Q_1(\nu)Q_2(\nu)\\
&+\tau_2\partial_{2}(\partial_{1}(Q_2(\nu))){}^{s_{2}}Q_1(\nu)Q_2(\nu)+\tau_2\tau_1{}^{s_{1}}Q_2(\nu)\partial_{2}(Q_1(\nu))Q_2(\nu)\\
&+\tau_2\partial_{1}(Q_2(\nu))\partial_{2}(Q_1(\nu))Q_2(\nu),
\end{split}
\end{align*}
and \eqref{ERPQ} equals
\begin{align*}
\begin{split}
&-1+\tau_1Q_1(\nu)+\tau_2Q_2(\nu)-\tau_1Q_1(\nu)\tau_2Q_2(\nu)+\tau_1Q_1(\nu)-\tau_1Q_1(\nu)\tau_1Q_1(\nu)\\&-\tau_2Q_2(\nu)\tau_1Q_1(\nu)
+\tau_1Q_1(\nu)\tau_2Q_2(\nu)\tau_1Q_1(\nu)\\
=&-1+\tau_1Q_1(\nu)+\tau_2Q_2(\nu)-\tau_1\tau_2{}^{s_{2}}Q_1(\nu)Q_2(\nu)
-\tau_1\partial_{2}(Q_1(\nu))Q_2(\nu)+\tau_1Q_1(\nu)\\
&-\tau_1\partial_{1}(Q_1(\nu))Q_1(\nu)-\tau_2\tau_1{}^{s_{1}}Q_2(\nu)Q_1(\nu)
-\tau_2\partial_{1}(Q_2(\nu))Q_1(\nu)\\
&+\tau_1\tau_2\tau_1{}^{s_{1}s_{2}}Q_{1}(\nu){}^{s_{1}}Q_2(\nu)Q_1(\nu)+
\tau_1\tau_2\partial_{1}({}^{s_{2}}Q_1(\nu)){}^{s_{1}}Q_2(\nu)Q_1(\nu)\\
&+\tau_1\partial_{1}(\partial_{2}(Q_1(\nu))){}^{s_{1}}Q_2(\nu)Q_1(\nu)+\tau_1\tau_2{}^{s_{2}}Q_1(\nu)\partial_{1}(Q_2(\nu))Q_1(\nu)\\
&+\tau_1\partial_{2}(Q_1(\nu))\partial_{1}(Q_2(\nu))Q_1(\nu).
\end{split}
\end{align*}

Since $r=s=t$, we have $\tau_2\tau_1\tau_2=\tau_1\tau_2\tau_1,$ and their corresponding coefficients are equal to each other
by \eqref{Q-relations1}. For the $\tau_1$-coefficients, it suffices to prove that
\begin{equation}
1-\partial_{1}(Q_1(\nu))+\partial_{1}(\partial_{2}(Q_1(\nu))){}^{s_{1}}Q_2(\nu)+\partial_{2}(Q_1(\nu))\partial_{1}(Q_2(\nu))=0.\label{rst-tau1}
\end{equation}
By an explicit calculation, we have $\partial_{1}(Q_1(\nu))=q+1,$ $\partial_{1}(\partial_{2}(Q_1(\nu)))=0,$ $\partial_{2}(Q_1(\nu))=-q,$ $\partial_{1}(Q_2(\nu))=-1.$ Thus, \eqref{rst-tau1} holds. The $\tau_2$-coefficients can be handed similarly. All the other coefficients are routine.

If $j_1, j_2, j_3$ are not all equal, by definition, we have the following identity:
\begin{align*}
g_{i}\diamond 1_{(\nu, j)}=\left\{
\begin{array}{ll}
(\tau_{i}Q_{i}(\nu, j)-P_{i}(\nu, j))1_{(\nu, j)}&\text{if $j_{i}=j_{i+1}$},\\[0.3em]
f_{i, j}\tau_{i}1_{(\nu, j)}&\text{otherwise.}
\end{array}
\right.
\end{align*}
For simplicity, we write $(i_1, i_2, i_3)$ instead of $((\nu_{i_1}, \nu_{i_2}, \nu_{i_3}),(j_{i_1}, j_{i_2}, j_{i_3}))$ or $(j_{i_1}, j_{i_2}, j_{i_3}).$ We need to check that $g_{2}g_{1}g_{2}\diamond 1_{(\nu, j)}=g_{1}g_{2}g_{1}\diamond 1_{(\nu, j)}.$ Then there are four cases for consideration.

\subsubsection*{Case 1: $j_1=j_2\neq j_3$.}

\

\vspace{0.35cm}\noindent By \eqref{rel-def-QHA2}, we have
\begin{align*}
\begin{split}
g_{2}g_{1}g_{2}\diamond 1_{(\nu, j)}&=g_{2}g_{1}\diamond f_{2, (123)}\tau_{2}1_{(123)}\\
&=f_{2, (123)}g_{2}\diamond f_{1, (132)}\tau_{1}\tau_{2}1_{(123)}\\
&=f_{2, (123)}f_{1, (132)}(\tau_{2}Q_{2}(312)-P_{2}(312))\tau_{1}\tau_{2}1_{(123)}\\
&=f_{2, (123)}f_{1, (132)}(\tau_{2}\tau_{1}\tau_{2}{}^{s_{2}s_{1}}Q_{2}(312)-\tau_{1}\tau_{2}{}^{s_{2}s_{1}}P_{2}(312))1_{(123)},
\end{split}
\end{align*}
and similarly, we have
\begin{align*}
\begin{split}
g_{1}g_{2}g_{1}\diamond 1_{(\nu, j)}=f_{2, (213)}f_{1, (231)}\tau_{1}\tau_{2}\tau_{1}Q_{1}(123)1_{(123)}-f_{2, (123)}f_{1, (132)}\tau_{1}\tau_{2}P_{1}(123))1_{(123)}.
\end{split}
\end{align*}
By the assumption on $j_1, j_2$ and $j_3$, we have $\tau_{2}\tau_{1}\tau_{2}=\tau_{1}\tau_{2}\tau_{1}.$ Using the fact that $f_{i, j}=f_{i+1, s_{i}s_{i+1}j}$ and \eqref{P-relations2} and \eqref{Q-relations1}, we see that the corresponding $\tau_{2}\tau_{1}\tau_{2}$-coefficients and $\tau_{1}\tau_{2}$-coefficients are equal to each other in the above two expressions.

\subsubsection*{Case 2: $j_1\neq j_2=j_3$.}

\

\vspace{0.3cm}\noindent Similar.

\subsubsection*{Case 3: $j_1=j_3\neq j_2$.}

\

\vspace{0.3cm}\noindent By the hypothesis, we have $\tau_{1}^{2}1_{(123)}=\tau_{2}^{2}1_{(123)}=1_{(123)}$ and $\tau_{2}\tau_{1}\tau_{2}=\tau_{1}\tau_{2}\tau_{1}.$ Then this case can also be checked by a routine calculation as in \emph{Case 1}.

\subsubsection*{Case 4: $j_1, j_2, j_3$ are all different.}

\

\vspace{0.3cm}\noindent This case can be easily checked.

Next we verify that the actions under $\Psi$ of the generators of $H_{n}(\Gamma)$ satisfy the relations \eqref{rel-def-QHA1}-\eqref{rel-def-QHA6}.

We can see that \eqref{rel-def-QHA1} holds by definition.

The relation \eqref{rel-def-QHA2} easily follows from \eqref{xyyx} and \eqref{phie-ephi}.

The relation \eqref{rel-def-QHA3} easily follows from \eqref{xyyx}, \eqref{phi-phi} and \eqref{phi-q}.

Next we check the relation \eqref{rel-def-QHA4}. If $a=i+1$, we need to prove that
\begin{align*}
\tau_{i}x_{i+1}1_{(\nu, j)}\ast n-x_{i}\tau_{i}1_{(\nu, j)}\ast n=\left\{
\begin{array}{ll}
1_{(\nu, j)}\ast n&\text{if $(\nu_{i}, j_{i})=(\nu_{i+1}, j_{i+1})$},\\[0.3em]
0&\text{otherwise,}
\end{array}
\right.
\end{align*}
for all $n\in N$ with an $N\in \mathcal{C}_{K}.$ To simplify notation for the remainder of the proof, we no longer write the elements ``$\ast n$" on the right hand side of all expressions, but remember it is always there.

If $(\nu_{i}, j_{i})=(\nu_{i+1}, j_{i+1}),$ we have
\begin{align*}
\tau_{i}x_{i+1}1_{(\nu, j)}=&(g_{i}+1)Q_{i}'(\nu, j)^{-1}(\nu_{i+1}^{-1}X_{i+1}-1)e(\nu, j)\\
=&((\nu_{i+1}^{-1}X_{i}-1)(g_{i}+1)+\nu_{i}^{-1}(qX_{i+1}-X_{i}))Q_{i}'(\nu, j)^{-1}e(\nu, j)\\
=&(x_{i}\tau_{i}+1)1_{(\nu, j)}.
\end{align*}

If $j_{i}=j_{i+1}$ and $\nu_{i}\neq \nu_{i+1}$, we have
\begin{align*}
\tau_{i}x_{i+1}1_{(\nu, j)}=&\Phi_{i}Q_{i}'(\nu, j)^{-1}(\nu_{i+1}^{-1}X_{i+1}-1)e(\nu, j)\\
=&(\nu_{i+1}^{-1}X_{i}-1)\Phi_{i}Q_{i}'(\nu, j)^{-1}e(\nu, j)\text{ by }\eqref{phi-x-nuj1}\\
=&x_{i}\tau_{i}1_{(\nu, j)}.
\end{align*}

If $j_{i}\neq j_{i+1}$, we have
\begin{align*}
\tau_{i}x_{i+1}1_{(\nu, j)}=&g_{i}\cdot f_{i, j}^{-1}(\nu_{i+1}^{-1}X_{i+1}-1)e(\nu, j)\\
=&(\nu_{i+1}^{-1}X_{i}-1)g_{i}\cdot f_{i, j}^{-1}e(\nu, j)\text{ by }\eqref{gxe-xge}\\
=&x_{i}\tau_{i}1_{(\nu, j)}.
\end{align*}

For $a=i,$ the proof is similar.

If $a\neq i, i+1,$ we have $\tau_{i}x_{a}1_{(\nu, j)}-x_{s_{i}a}\tau_{i}1_{(\nu, j)}=0$ by \eqref{phi-x}.

Since we have verified \eqref{rel-def-QHA4}, we can conclude that the identity \eqref{commu-relations} holds. For the relation \eqref{rel-def-QHA5}, by \eqref{rel-def-QHA2}, we have
\begin{align*}
\tau_{i}^{2}1_{(\nu, j)}=\Phi_{i}Q_{i}'(s_{i}\cdot(\nu, j))^{-1}(\tau_{i}1_{(\nu, j)}).
\end{align*}

If $(\nu_{i}, j_{i})=(\nu_{i+1}, j_{i+1}),$ we have
\begin{align*}
\tau_{i}^{2}1_{(\nu, j)}=&(g_{i}+1)(q\nu_{i+1}^{-1}X_{i+1}-\nu_{i}^{-1}X_{i})^{-1}(\tau_{i}1_{(\nu, j)})\\
=&(g_{i}+1)((q(x_{i+1}+1)-(x_{i}+1))^{-1}\tau_{i}1_{(\nu, j)}).
\end{align*}
Using \eqref{commu-relations}, we get
\begin{align*}
(q(x_{i+1}+1)-(x_{i}+1))^{-1}\tau_{i}1_{(\nu, j)}=&\tau_{i}(q(x_{i}+1)-(x_{i+1}+1))^{-1}1_{(\nu, j)}
-(q+1)\\\times(q(x_{i+1}+&1)-(x_{i}+1))^{-1}(q(x_{i}+1)-(x_{i+1}+1))^{-1}1_{(\nu, j)}.
\end{align*}
Thus, we have
\begin{align*}
\tau_{i}^{2}1_{(\nu, j)}=&(g_{i}+1)(\tau_{i}(q(x_{i}+1)-(x_{i+1}+1))^{-1}1_{(\nu, j)}
-(q+1)\\&\times(q(x_{i+1}+1)-(x_{i}+1))^{-1}(q(x_{i}+1)-(x_{i+1}+1))^{-1}1_{(\nu, j)})\\
=&(g_{i}+1)((g_{i}+1)Q_{i}'(\nu, j)^{-1}(q\nu_{i}^{-1}X_{i}-\nu_{i+1}^{-1}X_{i+1})^{-1}-(q+1)\\
&\times(q\nu_{i+1}^{-1}X_{i+1}-
\nu_{i}^{-1}X_{i})^{-1}(q\nu_{i}^{-1}X_{i}-\nu_{i+1}^{-1}X_{i+1})^{-1})e(\nu, j)\\
=&(g_{i}+1)(g_{i}-q)Q_{i}'(\nu, j)^{-1}(q\nu_{i}^{-1}X_{i}-\nu_{i+1}^{-1}X_{i+1})^{-1}e(\nu, j)\\
=&0.
\end{align*}

If $j_{i}=j_{i+1}$ and $\nu_{i}\neq \nu_{i+1}$, by \eqref{commu-relations}, we have
\begin{align*}
\tau_{i}^{2}1_{(\nu, j)}=&\Phi_{i}(Q_{i}(s_{i}\cdot(\nu, j))^{-1}\tau_{i}1_{(\nu, j)})\\
=&\Phi_{i}(\tau_{i}{}^{s_{i}}Q_{i}(s_{i}\cdot(\nu, j))^{-1}1_{(\nu, j)})\\
=&\Phi_{i}^{2}Q_{i}'(\nu, j)^{-1}{}^{s_{i}}Q_{i}'(s_{i}\cdot(\nu, j))^{-1}e(\nu, j).
\end{align*}
Since
\begin{align}
(1-P_{i}'(\nu, j))(q+P_{i}'(\nu, j))e(\nu, j)=&(1-\frac{1-q}{1-X_{i}X_{i+1}^{-1}})(q+\frac{1-q}{1-X_{i}X_{i+1}^{-1}})e(\nu, j)\notag\\
=&\frac{(qX_{i+1}-X_{i})(X_{i+1}-qX_i)}{(X_{i+1}-X_{i})^{2}}e(\nu, j)\label{p-ralation}\\
=&\Phi_{i}^{2}e(\nu, j)\text{ by }\eqref{phi-2-nuj1}\notag,
\end{align}
we see that \eqref{rel-def-QHA5} holds by \eqref{Q-relations2}.

If $j_{i}\neq j_{i+1}$, we have
\begin{align*}
\tau_{i}^{2}1_{(\nu, j)}=&g_{i}\cdot f_{i, s_{i}(j)}^{-1}(\tau_{i}1_{(\nu, j)})\\
=&f_{i, s_{i}(j)}^{-1}g_{i}\cdot g_{i}\cdot f_{i, j}^{-1}e(\nu, j)\\
=&(f_{i, j}f_{i, s_{i}(j)})^{-1}g_{i}^{2}e(\nu, j)\\
=&q^{-1}\cdot qe(\nu, j)=1_{(\nu, j)}.
\end{align*}

Finally we prove \eqref{rel-def-QHA6}. Without loss of generality we assume that $i=1$ and $n=3,$ and we need to show that \eqref{rel-def-QHA6} holds for any $\alpha=(\nu, j)\in K^{3}$ with $\nu=(\nu_1, \nu_2, \nu_3)$ and $j=(j_1, j_2, j_3).$

If $j_1=j_2=j_3$, we set $i:=\nu_1, l:=\nu_2, k:=\nu_3.$ We stop writing $j=(j_1, j_2, j_3)$ in the expressions of $1_{(\nu, j)}$, $e(\nu, j)$, $Q_{a}(\nu, j)$ and $P_{a}(\nu, j)$, that is, we write $1_{\nu}$ instead of $1_{(\nu, j)}$, $e(\nu)$ instead of $e(\nu, j),$ $Q_{a}(\nu)$ instead of $Q_{a}(\nu, j)$ and $P_{a}(\nu)$ instead of $P_{a}(\nu, j)$ for $a=1,2.$ To show that \eqref{rel-def-QHA6} holds, it suffices to consider the following five cases.

\subsubsection*{Case 1: $i, l, k$ are all different.}
\subsubsection*{Case 2: $i=l\neq k$.}
\subsubsection*{Case 3: $i\neq l=k$.}

\

\vspace{0.3cm}\noindent In these three cases, we need to show that $\tau_{1}\tau_{2}\tau_{1}1_{\nu}=\tau_{2}\tau_{1}\tau_{2}1_{\nu},$ which can be proved in exactly the same way as in the proof of [BK1, Theorem 4.2, p. 475] by using \eqref{phi-braid}. We omit the details.

\subsubsection*{Case 4: $i=k\neq l$.}

\

\vspace{0.3cm}\noindent In this case, we have
\begin{align*}
&\tau_{1}\tau_{2}\tau_{1}1_{\nu}\\
=&\Phi_{1}\Phi_{2}\Phi_{1}((Q_{1}(ili){}^{s_{1}}Q_{2}(lii){}^{s_{1}s_{2}}Q_{1}(lii))^{-1}1_{\nu})
+\Phi_{1}^{2}(Q_{1}(ili)^{-1}({}^{s_{1}}\partial_{2}(Q_1(lii)^{-1}))1_{\nu}),
\end{align*}
and
\begin{align*}
&\tau_{2}\tau_{1}\tau_{2}1_{\nu}\\
=&\Phi_{2}\Phi_{1}\Phi_{2}((Q_{2}(ili){}^{s_{2}}Q_{1}(iil){}^{s_{2}s_{1}}Q_{2}(iil))^{-1}1_{\nu})
+\Phi_{2}^{2}(Q_{2}(ili)^{-1}({}^{s_{2}}\partial_{1}(Q_2(iil)^{-1}))1_{\nu}).
\end{align*}
Using \eqref{Q-relations1}, \eqref{phi-2-nuj1} and \eqref{phi-braid}, we get that $\tau_{1}\tau_{2}\tau_{1}1_{\nu}-\tau_{2}\tau_{1}\tau_{2}1_{\nu}=A+B-C$, where
\begin{align*}
A=&(1-q)^2\frac{(X_1X_{3}-X_{2}^2)(X_1 X_{2}-qX_{2}X_{3})}{(X_1-X_{2})^2(X_{2}-X_{3})^2}((Q_{1}(ili){}^{s_{1}}Q_{2}(lii)Q_{2}(ili))^{-1}1_{\nu}),\\
B=&\frac{(X_{2}-qX_1)(X_1-qX_{2})}{(X_{2}-X_{1})(X_{1}-X_{2})}(Q_{1}(ili)^{-1}({}^{s_{1}}\partial_{2}(Q_1(lii)^{-1}))1_{\nu}),\\
C=&\frac{(X_{3}-qX_2)(X_2-qX_{3})}{(X_{3}-X_{2})(X_{2}-X_{3})}(Q_{2}(ili)^{-1}({}^{s_{2}}\partial_{1}(Q_2(iil)^{-1}))1_{\nu}).
\end{align*}
Noting that $Q_{2}(lii)1_{\nu}=i^{-1}(qX_{3}-X_{2})e(\nu),$ we get that
\begin{align*}
A=-i(1-q)^2\frac{(X_1X_{3}-X_{2}^2)X_{2}}{(X_1-X_{2})^2(X_{2}-X_{3})^2}((Q_{1}(ili)Q_{2}(ili))^{-1}1_{\nu}).
\end{align*}
Since
\begin{align*}\partial_{2}(Q_1(lii)^{-1})1_{\nu}=\frac{Q_1(lii)^{-1}-{}^{s_{2}}Q_1(lii)^{-1}}{x_{3}-x_{2}}1_{\nu}=
\frac{Q_1(lii)^{-1}-{}^{s_{1}}Q_2(ili)^{-1}}{x_{3}-x_{2}}1_{\nu},
\end{align*}
we have
\begin{align*}
B=\frac{(X_{2}-qX_1)(X_1-qX_{2})}{(X_{2}-X_{1})(X_{1}-X_{2})}\big(\frac{Q_{1}(ili)^{-1}Q_{2}(ili)^{-1}-
Q_{1}(ili)^{-1}({}^{s_{1}}Q_{1}(lii)^{-1})}{x_{1}-x_{3}}1_{\nu}\big).
\end{align*}
Similarly, we have
\begin{align*}
C=\frac{(X_{3}-qX_2)(X_2-qX_{3})}{(X_{3}-X_{2})(X_{2}-X_{3})}\big(\frac{Q_{1}(ili)^{-1}Q_{2}(ili)^{-1}-
Q_{2}(ili)^{-1}({}^{s_{2}}Q_{2}(iil)^{-1})}{x_{1}-x_{3}}1_{\nu}\big).
\end{align*}
Note that $(x_{1}-x_{3})1_{\nu}=i^{-1}(X_{1}-X_{3})e(\nu).$ By a direct computation we have the identity:
\begin{align*}
-i(1-q)^2& \frac{(X_1X_3-X_2^2)X_2}{(X_1-X_2)^2
(X_2-X_3)^2}
-\frac{i(X_2-qX_1)(X_1-qX_2)}{(X_1-X_2)^2(X_1-X_3)}\\
+&\frac{i(X_3-qX_2)(X_2-qX_3)}{(X_2-X_3)^2(X_1-X_3)}=0.
\end{align*}
Thus, we get that
\begin{align*}
A+B-C=&-\frac{(X_{2}-qX_1)(X_1-qX_{2})}{(X_{2}-X_{1})(X_{1}-X_{2})}\big(\frac{(Q_{1}(ili){}^{s_{1}}Q_{1}(s_{1}\cdot(ili)))^{-1}}{x_{1}-x_{3}}1_{\nu}\big)\\
&+\frac{(X_{3}-qX_2)(X_2-qX_{3})}{(X_{3}-X_{2})(X_{2}-X_{3})}\big(\frac{(Q_{2}(ili){}^{s_{2}}Q_{2}(s_{2}\cdot(ili)))^{-1}}{x_{1}-x_{3}}1_{\nu}\big).
\end{align*}
By using \eqref{Q-relations2} and \eqref{p-ralation}, we deduce that this equals 0 if $i\nslash l$, $1_{\nu}$ if $i\rightarrow l$, $-1_{\nu}$ if $i\leftarrow l$, and $(-2x_2+x_1+x_3)1_{\nu}$ if $i\rightleftarrows l.$ Thus, \eqref{rel-def-QHA6} holds.

\subsubsection*{Case 5: $i=l=k$.}

\

\vspace{0.3cm}\noindent This case is left as an exercise in the proof of [BK1, Theorem 4.2, p. 476]. Here we shall give a brief proof. In this case, by \eqref{Q-relations1}, we have
\begin{align*}
\tau_{1}\tau_{2}\tau_{1}1_{\nu}=&\Phi_{1}\Phi_{2}\Phi_{1}((Q_{1}(iii){}^{s_{2}}Q_{1}(iii)Q_{2}(iii))^{-1}1_{\nu})\\
&+\Phi_{1}^{2}(Q_{1}(iii)^{-1}({}^{s_{1}}\partial_{2}(Q_1(iii)^{-1}))1_{\nu})+\Phi_{1}(\partial_{1}(\partial_{2}(Q_1(iii)^{-1}))1_{\nu}),
\end{align*}
and
\begin{align*}
\tau_{2}\tau_{1}\tau_{2}1_{\nu}=&\Phi_{2}\Phi_{1}\Phi_{2}((Q_{2}(iii){}^{s_{1}}Q_{2}(iii)Q_{1}(iii))^{-1}1_{\nu})\\
&+\Phi_{2}^{2}(Q_{2}(iii)^{-1}({}^{s_{2}}\partial_{1}(Q_2(iii)^{-1}))1_{\nu})+\Phi_{2}(\partial_{2}(\partial_{1}(Q_2(iii)^{-1}))1_{\nu}).
\end{align*}
Using \eqref{Q-relations1}, \eqref{phi-braid} and the fact that $\Phi_{1}^{2}e(\nu)=(1+q)\Phi_{1}e(\nu)$ and $\Phi_{2}^{2}e(\nu)=(1+q)\Phi_{2}e(\nu),$ we see that $\tau_{1}\tau_{2}\tau_{1}1_{\nu}-\tau_{2}\tau_{1}\tau_{2}1_{\nu}=\Phi_{1}(D1_{\nu})-\Phi_{2}(E1_{\nu}),$
where
\begin{align*}
D=q(Q_{1}(iii){}^{s_{2}}Q_{1}(iii)Q_{2}(iii))^{-1}+(1+q)Q_{1}(iii)^{-1}({}^{s_{1}}\partial_{2}(Q_1(iii)^{-1}))
+\partial_{1}(\partial_{2}(Q_1(iii)^{-1})),
\end{align*}
and
\begin{align*}
E=q(Q_{1}(iii){}^{s_{2}}Q_{1}(iii)Q_{2}(iii))^{-1}+(1+q)Q_{2}(iii)^{-1}({}^{s_{2}}\partial_{1}(Q_2(iii)^{-1}))
+\partial_{2}(\partial_{1}(Q_2(iii)^{-1})).
\end{align*}
In order to show that $\tau_{1}\tau_{2}\tau_{1}1_{\nu}=\tau_{2}\tau_{1}\tau_{2}1_{\nu},$ it suffices to prove that $D=E=0.$ By a direct computation, we have
\begin{align*}
q(Q_{1}(iii){}^{s_{2}}Q_{1}(iii)Q_{2}(iii))^{-1}=\frac{q}{(qx_2-x_1)(qx_3-x_2)(qx_3-x_1)},
\end{align*}
\begin{align*}
(1+q)Q_{1}(iii)^{-1}({}^{s_{1}}\partial_{2}(Q_1(iii)^{-1}))=\frac{q(1+q)}{(qx_2-x_1)(qx_3-x_2)(qx_1-x_2)},
\end{align*}
\begin{align*}
\partial_{1}(\partial_{2}(Q_1(iii)^{-1}))=\frac{q(x_1+x_2-qx_3-q^{2}x_3)}{(qx_2-x_1)(qx_3-x_1)(qx_1-x_2)(qx_3-x_2)}.
\end{align*}
Then it is easy to see that $D=0.$ Similarly, we have $E=0.$ Thus, \eqref{rel-def-QHA6} holds.

If $j_1, j_2, j_3$ are not all equal, by definition, we have the following identity:\vspace{0.1cm}
\begin{align*}
\tau_{i}1_{(\nu, j)}=\left\{
\begin{array}{ll}
(g_{i}+(1-q)(1-X_{i}X_{i+1}^{-1})^{-1})Q_{i}'(\nu, j)^{-1}e(\nu, j)&\text{if $\nu_{i}\neq \nu_{i+1}$ and $j_{i}=j_{i+1}$},\\[0.3em]
(g_{i}+1)Q_{i}'(\nu, j)^{-1}e(\nu, j)&\text{if $\nu_{i}=\nu_{i+1}$ and $j_{i}=j_{i+1}$},\\[0.3em]
f_{i, j}^{-1}g_{i}e(\nu, j)&\text{if $j_{i}\neq j_{i+1}$.}
\end{array}
\right.
\end{align*}\vspace{0.2cm}
For simplicity, we write $(i_1, i_2, i_3)$ instead of $((\nu_{i_1}, \nu_{i_2}, \nu_{i_3}),(j_{i_1}, j_{i_2}, j_{i_3}))$ or $(j_{i_1}, j_{i_2}, j_{i_3}).$ We need to check that $\tau_{1}\tau_{2}\tau_{1}1_{(\nu, j)}=\tau_{2}\tau_{1}\tau_{2}1_{(\nu, j)}.$ Then there are four cases for consideration.
\subsubsection*{Case 1: $j_1=j_2\neq j_3$.}

\

\vspace{0.35cm}\noindent  In this case, we need to consider the following two subcases.

$\clubsuit$ If $\nu_1\neq \nu_2,$ we have
\begin{align*}
&\tau_{2}\tau_{1}\tau_{2}1_{(123)}\\
=&((g_{2}+(1-q)(1-X_{2}X_{3}^{-1})^{-1})Q_{2}'(312)^{-1}e(312))(\tau_{1}\tau_{2}1_{(123)})\\
=&(g_{2}+(1-q)(1-X_{2}X_{3}^{-1})^{-1})(Q_{2}(312)^{-1}\tau_{1}\tau_{2}1_{(123)})\\
=&(g_{2}+(1-q)(1-X_{2}X_{3}^{-1})^{-1})(\tau_{1}\tau_{2}{}^{s_{2}s_{1}}Q_{2}(312)^{-1}1_{(123)})\\
=&(g_{2}+(1-q)(1-X_{2}X_{3}^{-1})^{-1})(\tau_{1}1_{(132)})(\tau_{2}Q_{1}(123)^{-1}1_{(123)})\\
=&(g_{2}+(1-q)(1-X_{2}X_{3}^{-1})^{-1})(f_{1, (132)}^{-1}g_{1}e(132))(f_{2, (123)}^{-1}g_{2}e(123))(Q_{1}(123)^{-1}1_{(123)})\\
=&f_{1, (132)}^{-1}f_{2, (123)}^{-1}(g_{2}g_{1}g_{2}+(1-q)(1-X_{2}X_{3}^{-1})^{-1}g_{1}g_{2})(Q_{1}(123)^{-1}1_{(123)})\\
=&f_{1, (132)}^{-1}f_{2, (123)}^{-1}(g_{1}g_{2}g_{1}+(1-q)g_{1}g_{2}(1-X_{1}X_{2}^{-1})^{-1})(Q_{1}(123)^{-1}1_{(123)})\text{ by }\eqref{gxe-xge}\\
=&f_{1, (132)}^{-1}f_{2, (123)}^{-1}g_{1}g_{2}(g_{1}+(1-q)(1-X_{1}X_{2}^{-1})^{-1})Q_{1}'(123)^{-1}e(123),
\end{align*}
and
\begin{align*}
\tau_{1}\tau_{2}\tau_{1}1_{(123)}=&(\tau_{1}1_{(231)})(\tau_{2}1_{(213)})(\tau_{1}1_{(123)})\\
=&f_{1, (231)}^{-1}f_{2, (213)}^{-1}g_{1}g_{2}(g_{1}+(1-q)(1-X_{1}X_{2}^{-1})^{-1})Q_{1}'(123)^{-1}e(123).
\end{align*}
By using the fact that $f_{i, j}=f_{i+1, s_{i}s_{i+1}j}$, we can get that $\tau_{1}\tau_{2}\tau_{1}1_{(123)}=\tau_{2}\tau_{1}\tau_{2}1_{(123)}$.

$\clubsuit$ If $\nu_1=\nu_2,$ we have
\begin{align*}
\tau_{2}1_{(312)}=(g_{2}+1)Q_{2}'(312)^{-1}e(312),
\end{align*}
\begin{align*}
\tau_{1}1_{(123)}=(g_{1}+1)Q_{1}'(123)^{-1}e(123).
\end{align*}
We can get the desired identity by the same calculation as above.
\subsubsection*{Case 2: $j_1\neq j_2=j_3$.}

\

\vspace{0.3cm}\noindent Similar.
\subsubsection*{Case 3: $j_1=j_3\neq j_2$.}

\

\vspace{0.3cm}\noindent In this case, we also need to consider the following two subcases: $\nu_1\neq \nu_2$ and $\nu_1=\nu_2.$ Note that $g_{1}^{2}e(123)=g_{2}^{2}e(123)=qe(123).$ Then this case can also be checked by a routine calculation as in \emph{Case 1}.

\subsubsection*{Case 4: $j_1, j_2, j_3$ are all different.}

\

\vspace{0.3cm}\noindent This case can be easily checked.

It is obvious that $\Phi\circ\Psi=$Id and $\Psi\circ\Phi=$Id. Thus, $\Phi$ and $\Psi$ establish an equivalence of categories.

\section{The degenerate case}

In this section, we consider the case $q=1.$ We introduce a degenerate affine Yokonuma-Hecke algebra, and establish an equivalence between its module category and its suitable counterpart for the quiver Hecke algebra. Since most of the calculations are entirely similar to the non-degenerate case, we shall not write them explicitly, and just state the main result.

Following [Ro, $\S5.1$], we define the degenerate affine Yokonuma-Hecke algebra $\widehat{Y}_{r,n}(1).$

\begin{definition}
The degenerate affine Yokonuma-Hecke algebra, denoted by $\widehat{Y}_{r,n}(1)$, is an $\mathcal{R}$-associative algebra generated by the elements $t_{1},\ldots,t_{n},f_{1},\ldots,f_{n-1},y_1,\ldots, y_{n}$ in which the generators $t_{1},\ldots,t_{n},g_{1},$ $\ldots,g_{n-1}$ satisfy the following relations:
\begin{align}
f_if_j&=f_jf_i \qquad \qquad\qquad\quad\hspace{0.3cm}\mbox{for all $i,j=1,\ldots,n-1$ such that $\vert i-j\vert \geq 2$;}\label{drel-def-Y1}\\[0.1em]
f_if_{i+1}f_i&=f_{i+1}f_if_{i+1} \qquad \quad\qquad\hspace{0.05cm}\mbox{for all $i=1,\ldots,n-2$;}\label{drel-def-Y2}\\[0.1em]
t_it_j&=t_jt_i \qquad\qquad\qquad\qquad  \mbox{for all $i,j=1,\ldots,n$;}\label{drel-def-Y3}\\[0.1em]
f_it_j&=t_{s_i(j)}f_i \quad \quad\qquad\qquad\hspace{0.25cm}\mbox{for all $i=1,\ldots,n-1$ and $j=1,\ldots,n$;}\label{drel-def-Y4}\\[0.1em]
t_i^r&=1 \quad \qquad\qquad\qquad\qquad\mbox{for all $i=1,\ldots,n$;}\label{drel-def-Y5}\\[0.2em]
f_{i}^{2}&=1 \qquad \hspace{2.71cm}\mbox{for all $i=1,\ldots,n-1$;}\label{drel-def-Y6}
\end{align}
together with the following relations concerning the generators $y_1,\ldots, y_{n}$:
\begin{align}
y_{i}y_{j}&=y_{j}y_{i};\label{drel-def-Y7}\\[0.1em]
f_{i}y_{i+1}&=y_{i}f_{i}+e_{i};\label{drel-def-Y8}\\[0.1em]
f_{i}y_{j}&=y_{j}f_{i} \qquad \qquad \quad\mbox{for all $j\neq i, i+1$;}\label{drel-def-Y9}\\[0.1em]
t_{j}y_{i}&=y_{i}t_{j} \hspace{0.07cm}\qquad \qquad\quad\mbox{for all $i, j=1,\ldots,n$;}\label{drel-def-Y10}
\end{align}
where for each $1\leq i\leq n-1$,
$$e_{i} :=\frac{1}{r}\sum\limits_{s=0}^{r-1}t_{i}^{s}t_{i+1}^{-s}.$$
\end{definition}

Let $\widehat{Y}_{r,n}^{k}(1)=k\otimes_{\mathcal{R}}\widehat{Y}_{r,n}(1)$ and let $M$ be a $\widehat{Y}_{r,n}^{k}(1)$-module. Given $\alpha=(\nu, j)\in K^{n},$ we denote by $M_{\alpha}$ the subspace of $M$ on which $y_{a}-\nu_{a}$ acts locally nilpotently for $1\leq a\leq n,$ and simultaneously, $t_{a}-\zeta^{j_{a}}$ acts as zero for $1\leq a\leq n.$ Let $\bar{\mathcal{C}}_{K}$ be the category of finitely generated $\widehat{Y}_{r,n}^{k}(1)$-modules $M$ such that $M=\bigoplus_{\alpha\in K^{n}}M_{\alpha}.$

Given an object $M\in \bar{\mathcal{C}}_{K}$, we can write $M$ as the direct sum of its weight spaces (simultaneous generalized eigenspaces):
\begin{equation}
M(\nu, j)=\{v\in M\:|\:(y_{a}-\nu_{a})^{N}v=(t_{a}-\zeta^{j_{a}})v=0\text{ for all }1\leq a\leq n\text{ and }N\gg 0\}.\label{M-J}
\end{equation}
Considering once again the weight space decomposition of the regular module, we deduce that there is a family $\{e(\nu,j)\}_{(\nu,j)\in K^{n}}$ of mutually orthogonal idempotents in $\widehat{Y}_{r,n}^{k}(1)$ such that $e(\nu,j)M=M(\nu, j)$ for each $M\in \bar{\mathcal{C}}_{K}$. In fact, each $e(\nu, j)$ lies in the commutative subalgebra generated by $y_{1},\ldots, y_{n}, t_{1},\ldots,t_{n}$, all but finitely many of the $e(\nu, j)$'s are zero, and their sum is the identity element in $\widehat{Y}_{r,n}^{k}(1).$

The proof of the next theorem is entirely similar to that of Theorem \ref{Morita-equi}, which requires extremely careful verification.
\begin{theorem}\label{dege-Morita-equi}
There is an equivalence of categories $\bar{\Phi}: \mathcal{C}_{\Gamma}^{0}\rightarrow \bar{\mathcal{C}}_{K}$ given by $M\mapsto M,$ and where for each $\alpha=(\nu, j)\in K^{n},$ $t_{a}$ acts on $1_{\alpha}M$ by $\zeta^{j_{a}},$ $y_{a}$ acts on $1_{\alpha}M$ by $x_{a}+\nu_{a}$ for $1\leq a\leq n,$ and $f_{i}$ acts on $1_{\alpha}M$ by $\tau_{i}q_{i}(\nu, j)-p_{i}(\nu, j)$ for $1\leq i\leq n-1,$ where

\[q_{i}(\nu, j)=
\begin{cases}
\begin{cases}
x_{i+1}-x_{i}+1 & \text{if } \nu_i = \nu_{i+1},
\\
\left(x_{i}-x_{i+1}-1\right)^{-1} & \text{if } \nu_{i+1} =\nu_{i}+1,
\\
\frac{x_{i+1}-x_{i}+\nu_{i+1}-\nu_{i}-1}
{x_{i}-x_{i+1}+\nu_{i}-\nu_{i+1}}
 & \text{otherwise},
\end{cases}
& \text{if } j_i = j_{i+1},
\\
1 & \text{if } j_i \neq j_{i+1},
\end{cases}
\]
\setlength{\parskip}{2pt}
and

\[p_{i}(\nu, j)=
\begin{cases}
\begin{cases}
1 & \text{if } \nu_i = \nu_{i+1},
\\
\frac{1}{x_{i}-x_{i+1}+\nu_{i}-\nu_{i+1}}
 & \text{otherwise},
\end{cases}
& \text{if } j_i = j_{i+1},
\\
0 & \text{if } j_i \neq j_{i+1}.
\end{cases}
\]

The inverse functor $\bar{\Psi}: \bar{\mathcal{C}}_{K}\rightarrow \mathcal{C}_{\Gamma}^{0}$ is given by $N\mapsto N,$ where for each $\alpha=(\nu, j)\in K^{n},$ $1_{\alpha}$ acts on $N$ by $e(\nu, j),$ $x_{a}$ acts on $N_{\alpha}$ by $y_{a}-\nu_{a}$ for $1\leq a\leq n,$ and $\tau_{i}$ acts on $N_{\alpha}$ by $\varphi_{i}q_{i}'(\nu, j)^{-1}=(f_{i}+p_{i}'(\nu, j))q_{i}'(\nu, j)^{-1}$ for $1\leq i\leq n-1,$ where $p_{i}'(\nu, j)$ and $q_{i}'(\nu, j)$ are defined by replacing $x_{i}$ with $y_{i}-\nu_{i}$ in the expressions of $p_{i}(\nu, j)$ and $q_{i}(\nu, j),$ and
\[\varphi_i=f_i + \sum_{\substack{(\nu, j) \in K^{n} \\ \nu_i \neq \nu_{i+1} \\ j_i = j_{i+1}}} \left(y_i-y_{i+1}\right)^{-1} e(\nu, j) + \sum_{\substack{(\nu, j) \in K^{n} \\ \nu_i = \nu_{i+1} \\ j_i = j_{i+1}}} e(\nu, j).
\]
\end{theorem}

\begin{remark}\label{Typos}
We point out that there are some typos in [Rou1, Theorems 3.16 and 3.19] (also [Rou2, Theorems 3.11 and 3.12]). For example, let us look at [Rou1, Theorem 3.19]. We claim that when $\nu_{i+1}=q\nu_{i},$ the formula is not right, since it does not satisfy
\begin{align}\label{T-2}
T_{i}^{2}\diamond 1_{\nu}=(q+(q-1)T_{i})\diamond 1_{\nu}
\end{align}
Set
\begin{align*}
Q_{i}(\nu)=
\begin{cases}
(q^{-1}x_{i}-x_{i+1})^{-1} & \text{if } \nu_{i+1}=q\nu_{i},
\\
\frac{q\nu_{i}x_{i}-\nu_{i+1}x_{i+1}}{\nu_{i}x_{i}-\nu_{i+1}x_{i+1}}
 & \text{if } \nu_{i}=q\nu_{i+1},
\end{cases}
\end{align*}
and
\begin{align*}
P_{i}(\nu)=(q-1)\nu_{i+1}x_{i+1}(\nu_{i}x_{i}-\nu_{i+1}x_{i+1})^{-1}\text{ if } \nu_{i+1}\neq \nu_{i}.
\end{align*}
When $\nu_{i+1}=q\nu_{i},$ if \eqref{T-2} holds, we deduce that
\begin{align*}
Q_{i}(s_{i}\cdot(\nu)){}^{s_{i}}Q_{i}(\nu)\tau_{i}^{2}1_{\nu}=(1-P_{i}(\nu))(q+P_{i}(\nu))1_{\nu}.
\end{align*}
Then we must have
\begin{align*}
(\nu_{i+1}&x_{i}-\nu_{i}x_{i+1})^{-1}(q\nu_{i+1}x_{i}-\nu_{i}x_{i+1})(q^{-1}x_{i+1}-x_{i})^{-1}(x_{i+1}-x_{i})\\
=&q+(1-q)(q-1)\nu_{i+1}x_{i+1}(\nu_{i}x_{i}-\nu_{i+1}x_{i+1})^{-1}-(q-1)^{2}\nu_{i+1}^{2}x_{i+1}^{2}(\nu_{i}x_{i}-\nu_{i+1}x_{i+1})^{-2}.
\end{align*}
Replacing $\nu_{i}$ with $q^{-1}\nu_{i+1},$ we get that
\begin{align*}
(q^{-1}x_{i+1}&-qx_{i})(x_{i+1}-x_{i})(q^{-1}x_{i+1}-x_{i})^{-2}\\
=&q-(q-1)^{2}x_{i+1}(q^{-1}x_{i}-x_{i+1})^{-1}-(q-1)^{2}x_{i+1}^{2}(q^{-1}x_{i}-x_{i+1})^{-2}.
\end{align*}
After a careful calculation, we deduce that $q+2q^{-2}=q^{-3}+2,$ which is a contradiction.

We also claim that when $\nu_{i+1}=\nu_{i},$ the action of $X_{i}$ given in [Rou1, Theorem 3.19] or [Rou2, Theorem 3.11] is not right, since it does not satisfy
\begin{align*}
T_{i}X_{i}T_{i}\diamond 1_{\nu}=qX_{i+1}\diamond 1_{\nu},
\end{align*}
which is equivalent to
\begin{align}\label{T-X-T}
X_{i}T_{i}\diamond 1_{\nu}=(T_{i}+(1-q))X_{i+1}\diamond 1_{\nu}.
\end{align}

When $X_{i}$ acts as $\nu_{i}(x_{i}+1)$, from \eqref{T-X-T} we easily deduce that $-\nu_{i}(qx_{i+1}+1)=-q\nu_{i}(x_{i+1}+1),$ which is a contradiction.

When $X_{i}$ acts as $x_{i}+\nu_{i}$, from \eqref{T-X-T} we easily deduce that $-qx_{i+1}-\nu_{i}=-qx_{i+1}-q\nu_{i},$ which is a contradiction.

There exist some similar typos in [Rou1, Theorem 3.16].
\end{remark}



\end{document}